\documentclass[twoside]{article}
\usepackage{amsmath,amssymb,graphicx}
\usepackage[ngerman,english]{babel}
\usepackage{mathptmx}
\usepackage{graphics}
\usepackage{graphicx}
\usepackage{epstopdf}
\usepackage{grffile}
\usepackage{float}
\usepackage{floatflt}
\usepackage{epsfig, epsf, rotating, fancyhdr, afterpage, multicol}
\usepackage[hidelinks]{hyperref}
\usepackage{multirow}
\usepackage{listings}
\usepackage{booktabs}

\usepackage{pdfpages}

\usepackage{xcolor}
\definecolor{darkblue}{rgb}{0,0,0.5} 
\usepackage{transparent}
\usepackage{cite}
\usepackage{amsmath,amssymb,amsbsy,amsfonts,amsthm, latexsym,amsopn,amstext,amsxtra,euscript,amscd}
\usepackage{anysize}
\usepackage[margin=12pt,labelfont=bf, format=hang, justification =raggedright]{subcaption}%font=scriptsize
\usepackage[margin=12pt,labelfont=bf, labelsep=endash, format=hang, justification =raggedright]{caption}
\usepackage{titlesec}
\titleformat*{\section}{\large\bfseries}
\titleformat*{\subsection}{\normalsize\bfseries}
\titleformat*{\subsubsection}{\normalsize\bfseries}

\numberwithin{equation}{section}
\numberwithin{remark}{section}
\numberwithin{definition}{section}
\numberwithin{lemma}{section}
\numberwithin{property}{section}

\newcommand{\sgn}{\operatorname{sgn}}
\begin{document}

%%%%%%%%% TITLE of abstract %%%%%%%%%

\begin{center}
\textbf{\large A Hybrid Riemann Solver for Large Hyperbolic Systems of Conservation Laws}
\end{center}

%\medskip

%%%%%%%%%  AUTHORS %%%%%%%%%

\begin{center}
Birte Schmidtmann$^{a, *}$, Manuel Torrilhon$^a$\\
{\small{${}^a$ \emph{Center for Computational Engineering Science, RWTH Aachen University, Schinkelstr. 2, 52062 Aachen, Germany}}}
%\texttt{schmidtmann@mathcces.rwth-aachen.de}
\end{center}

%\keywords{Finite Volume Method, incomplete Riemann solvers, conservation laws, hyperbolic systems, Euler equations, ideal Magnetohydrodynamics.}
%novel scientific contribution of this work: This work presents a new family of HLL-type Riemann solvers. These solvers require the same input information as HLL and one additional flux evaluation and yield solutions with less dissipation, compared to classical solvers. This is especially interesting for large systems of conservation laws and is demonstrated for ideal MHD and Grad's 13 moment equations.

\bigskip

%%%%%%%%%  text of ABSTRACT to be included below  %%%%%%
\begin{abstract}
\noindent We are interested in the numerical solution of large systems of hyperbolic conservation laws or systems in which the characteristic decomposition is expensive to compute. Solving such equations using finite volumes or discontinuous Galerkin requires a numerical flux function which solves local Riemann problems at cell interfaces. There are various methods to express the numerical flux function. On the one end, there is the robust but very diffusive Lax-Friedrichs solver; on the other end the upwind Godunov solver which respects all resulting waves. The drawback of the latter method is the costly computation of the eigensystem. 

This work presents a family of simple first order Riemann solvers, named HLLX$\omega$, which avoid solving the eigensystem. The new method reproduces all waves of the system with less dissipation than other solvers with similar input and effort, such as HLL and FORCE. The family of Riemann solvers can be seen as an extension or generalization of the methods introduced by Degond et al. \cite{DegondPeyrardRussoVilledieu1999}. We only require the same number of input values as HLL, namely the globally fastest wave speeds in both directions, or an estimate of the speeds. Thus, the new family of Riemann solvers is particularly efficient for large systems of conservation laws when the spectral decomposition is expensive to compute or no explicit expression for the eigensystem is available.
\end{abstract}
	%
%==============================
\section{Introduction}\label{sec:introduction}
%==============================
	%
	In finite volume methods, integrating conservation laws over a control volume leads to a formulation which requires the evaluation of local Riemann problems at each cell interface. The initial states for these problems are typically given by the left and right adjacent cell values. Since these local Riemann problems have to be solved many times in order to find the numerical solution, the Riemann solver is a building block of the finite volume method. Over the last decades, many different Riemann solvers were developed, see e.g. \cite{ToroRiemannSolvers} for a broad overview. The main challenges are the need for computational efficiency and easy implementation, while at the same time, accurate results without artificial oscillations need to be obtained.

Riemann solvers can be classified into complete and incomplete schemes, depending on whether all present characteristic fields are considered in the model or not. According to this classification, upwind, Godunov's method and Roe's scheme, are complete \cite{LeVeque1992}. These schemes yield monotone results, however, an evaluation of the eigensystem of the flux Jacobian is needed. Especially for large systems this characteristic decomposition is expensive to compute and in some cases an analytic expression is not available at all. Nevertheless, using Roe's scheme, all waves are well-resolved and it typically yields the best resolution of the Riemann wave fan. In order to reduce computational cost while keeping high resolution, there have been many attempts to approximate the upwind scheme without solving the eigenvalue problem, see e.g. \cite{DegondPeyrardRussoVilledieu1999, Torrilhon2012, CastroGallardoMarquina2014} and references therein.
	
	In this article, we are interested in incomplete Riemann solvers. In comparison to complete solvers, they need less characteristic information and are easier to implement. However, they contain more dissipation and thus, yield lower resolution, especially of slow waves. Nevertheless, in many test cases, these Riemann solvers may be sufficient to obtain good results, especially if the system contains only fast waves. 
	\newline \newline
	The rest of the paper is structured as follows. In Sec.~\ref{sec:FV} we introduce the necessary notation for finite volume (FV) schemes and Riemann problems in general. Sec.~\ref{sec:riemannsolvers} reviews some well-known Riemann solvers. In Sec.~\ref{sec:hybridSolvers} we discuss some hybrid Riemann solvers, i.e. solvers which are constructed as weighted combinations of the ones presented in the previous section. In Sec.~\ref{sec:JacobianFree} we describe how to implement the described flux functions without knowing the flux Jacobian, i.e. a Jacobian-free implementation. Then, in Sec.~\ref{sec:HLLXomega} we present the new family of Riemann solvers and discusses construction and parameter choices. The numerical results of Sec.~\ref{sec:numericalresults} underline the excellent performance of the new solvers and finally, in Sec.~\ref{sec:conclusion} we draw some conclusions. Appendix \ref{sec:Appendix} gives detailed information on how to implement the new solvers and appendix \ref{sec:AppendixB} show more solution plots of the numerical test cases.
%
%=============================================
\section{Finite Volume Method} \label{sec:FV}
%=============================================
%
We consider a Cauchy Problem of the hyperbolic conservation law
	\begin{align}	
		\label{eq:conservation_law}
		\begin{cases}	
			\partial_t U(x,t)+\partial_x f(U(x,t))\hspace{-0.3cm} &= 0, \qquad\qquad \text{in}\; \mathbb{R}\times \mathbb{R}^+\\
			U(x,0)&=U_0(x),
		\end{cases}	
	\end{align}			
	in one space dimension, equipped with the initial condition $U_0(x)$, where $U  = (U_1, \ldots, U_N )^T: \mathbb{R}\times \mathbb{R}^+ \rightarrow \mathbb{R}^N$. The flux function $f: \mathbb{R}^N\rightarrow \mathbb{R}^N$, as well as the initial condition $U_0$ are supposed to be given. We are interested in hyperbolic systems, i.e. the Jacobian matrix $A(U) =\partial f/\partial U$ is diagonalizable and has $N$ real-valued eigenvalues $\forall U\in\mathbb{R}^N$.
	We consider a regular grid in space, with the positions of the cell centers denoted by $x_i,\;i\in\mathbb{Z}$ and with uniform space intervals of size $ \Delta x$. The grid cells are defined by $C_i = [x_{i-1/2},\; x_{i+1/2}]$, where $x_{i\pm j} = x_i \pm j\Delta x$.
	Finite volume methods aim at approximating the cell averages of the true solution of \eqref{eq:conservation_law} with high accuracy, see e.g. \cite{LeVeque2002}. The cell average of the true solution $U(\cdot, t^n)$ in cell $C_i$ at time $t^n$ is given by 
	\begin{align}
		\label{eq:exSol}
		\bar U_i(t^n) = \frac{1}{\Delta x} \int_{C_i} U(x,t^n) dx.
	\end{align}
	The goal is to find an update rule to advance the approximate cell averages from a given time $t^n$ to a new time $t^{n+1} = t^n+\Delta t$, such that the true cell averages are well-approximated. In addition, the approximate solution should not develop any (relevant) spurious oscillations.
Integrating Eq.~\eqref{eq:conservation_law} over the cell $C_i$ and dividing by $\Delta x$ yields
	\begin{align}
		\label{eq:approxConsLaw}
		\frac{d \bar U_i}{d t} = -\frac{1}{\Delta x} \left[ f(U(x_{i + 1/2}, t)) - f(U(x_{i - 1/2}, t))\right]
	\end{align}
	which is still exact. We now want to find an approximation of the solution $\bar u_i^n$ satisfying $\bar u_i^n \approx \bar U_i(t^n)$. The quality of the approximation $\bar{u}_i$ depends on the accurate approximation of the fluxes at the cell boundaries $f(U(x_{i \pm 1/2},t))$.
This is achieved by constructing a numerical flux function $\hat{f}(u,v)$ which is Lipschitz continuous and consistent with the true flux function $\hat{f}(U, U)=f(U)$, see e.g.~\cite{LeVeque2002}. The numerical flux at the right boundary of cell $C_i$ is then given by 
	\begin{subequations}
		\label{eq:approxFlux}
		\begin{align}
			\hat{f}_{i+1/2}=\hat{f}\left(U^{(-)}_{i+1/2},U^{(+)}_{i+1/2}\right).
		\end{align}
	\end{subequations}	
	It takes as input variables the left and right limiting values of the solution vector $U$ at the interface $i+1/2$. One can construct higher order schemes by inserting higher-order reconstructions of these interface values, for example by using weighted essentially non-oscillatory (WENO) schemes \cite{LiuOsherChan1994, JiangShu1996} or limiting methods, e.g. \cite{VanLeer1979, ColellaWoodward1984, Marquina1994, CadaTorrilhon2009, SchmidtmannSeiboldTorrilhon2015}. The focus of this paper is the development of new Riemann solvers and not the order of accuracy. A comment on the extension of the newly developed schemes to higher order in space and time can be found in Sec.~\ref{subsec:higherOrder}. For the development of the methods however, the input values of the flux function are simply the left and right cell mean values, $U^{(-)}_{i+1/2}\equiv \bar{u}_i$ and $ U^{(+)}_{i+1/2}\equiv \bar{u}_{i+1}$. For the sake of simplicity, these are also denoted by $U_L$ and $U_R$. In summary, the evolution of cell averages is given by
	\begin{align}
		\label{eq:approxConsLawApprox}
		\frac{d \bar u_i}{d t} = -\frac{1}{\Delta x} \left( \hat{f}_{i+1/2} - \hat{f}_{i-1/2}\right).
	\end{align}
	As mentioned above, the numerical flux function is the crucial point of the finite volume method and determines the order of accuracy. The numerical flux function can be written in the general form
	\begin{align}
		\label{eq:gerneralNumFlux}
		\hat{f}(U_L, U_R) = \frac{1}{2}\left(f(U_L) + f(U_R) \right) - \frac{1}{2} D(U_L, U_R)\; (U_R-U_L)
	\end{align}
	with the dissipation matrix $D$. This matrix depends on the left and right adjoint states and determines the form of the numerical flux function $\hat f$. 
	
	The discussion of different Riemann Solvers is easier in the quasi-linear formulation involving the flux Jacobian $A(U) =\partial f/\partial U$ and some information on its characteristics. If the flux Jacobian is not available, we assume there exists a Roe Matrix $\tilde A$ satisfying $\tilde A(U_L, U_R) (U_L-U_R) = F(U_L) - F(U_R)$ amongst other properties \cite{Roe1981}. From now on, we write $\tilde{A}$ even though the flux Jacobian $A$ might exist. 
	
	We introduce some notation which is useful for the following discussions on Riemann solvers. Let us define the maximal and minimal eigenvalues of $\tilde{A}(U)$ by
	\begin{subequations}
		\label{eq:minmaxeigenvalues}
		\begin{align}
			\lambda_\text{min}(U) &= \min\{\lambda\;\;|\; \lambda\;\text{is an eigenvalue of}\;\tilde{A}(U)\} \\			
			\lambda_\text{max}(U) &= \max\{\lambda\;\;|\; \lambda\;\text{is an eigenvalue of}\;\tilde{A}(U)\}.
 		\end{align}
	\end{subequations}	
	The eigenvalues of $\tilde{A}(U)$ are also called the characteristic speeds of the corresponding system. From Eq.~\eqref{eq:minmaxeigenvalues} it follows that the spectral radius at $U$ is given by 
	\begin{align}
		\bar\lambda(U) = \max\{|\lambda_\text{min}(U)|, |\lambda_\text{max}(U)|\}.
		\label{eq:spectralRadius}
	\end{align} 

	We would like to point out that for general non-linear systems of equations, discontinuities propagate with wave speeds which are not equal to the eigenvalues. However, for $\|U_L-U_R\|$ sufficiently small we can assume that discontinuities propagate at the constant characteristic speeds $\lambda_i$ \cite{LeVeque1992}. Therefore, from now on we will use the terms 'wave speed' and 'characteristic speed' interchangeably and we identify the globally fastest wave speed with the spectral radius $\bar\lambda$.	
	
	The dissipation matrix $D$ depends on the left and right states $U_L$ and $U_R$, respectively. Therefore, $D$ also depends on characteristics of the flux Jacobian. As mentioned above, the flux Jacobian can be diagonalized as $\tilde{A}(U)=T(U)\,\Lambda(U)\,T(U)^{-1}$, where $\Lambda(U)$ is the eigenvalue matrix $\Lambda(U)=\text{diag}(\lambda_1(U),\ldots,\lambda_N(U))$, $\lambda_1<\lambda_2<\ldots<\lambda_N$ and $T(U)$ the corresponding eigenvector matrix.
	
	Since the dissipation matrix $D$ is a function of the flux Jacobian $\tilde{A}$, it can be shown that	
	\begin{align}
		\frac{\Delta t}{\Delta x}D(\tilde{A}) = T\,\frac{\Delta t}{\Delta x}\,\text{diag}(\tilde d(\lambda_1), \ldots, \tilde d(\lambda_N))\,T^{-1} = T\,\text{diag}(d(\nu_1), \ldots, d(\nu_N))\,T^{-1},
		\label{eq:trafoDd}
	\end{align}
	holds true for all dissipation matrices discussed in this paper. This can easily be seen, since all dissipation matrices considered in this work are either of polynomial nature, given by the absolute value of the flux Jacobian, or linear combinations of polynomials and the absolute value function. More specifically, in this paper it holds that
	\begin{align*}
		D(\tilde{A}) = \alpha_0 \tilde{A}^0 + \alpha_1 \tilde{A}^1 + \alpha_2 \tilde{A}^2+ \alpha_3 |\tilde{A}|
	\end{align*}
	with some coefficients $\alpha_i\in\mathbb{R},\;i=1,\ldots,4$, which might be zero. It then follows from Eq.~\eqref{eq:trafoDd} that
	\begin{align*}
		\frac{\Delta t}{\Delta x}D(\tilde{A}) &= \frac{\Delta t}{\Delta x}\,T\left( \alpha_0 I + \alpha_1 \Lambda + \alpha_2 \Lambda^2+ \alpha_3 |\Lambda|\right) T^{-1}\\		
		&= T\,\frac{\Delta t}{\Delta x}\text{diag}\left(\alpha_0 + \alpha_1 \lambda_1 + \alpha_2 \lambda_1^2+ \alpha_3 |\lambda_1|,\ldots, \alpha_0 + \alpha_1 \lambda_N + \alpha_2 \lambda_N^2+ \alpha_3 |\lambda_N|\right) T^{-1}\\
		&= T\,\frac{\Delta t}{\Delta x}\,\text{diag}(\tilde d(\lambda_1), \ldots, \tilde d(\lambda_N))\,T^{-1} \\
		&= T\,\text{diag}\left(\alpha_0\frac{\Delta t}{\Delta x} + \alpha_1 \nu_1 + \alpha_2\frac{\Delta x}{\Delta t} \nu_1^2+ \alpha_3 |\nu_1|,\ldots, \alpha_0\frac{\Delta t}{\Delta x} + \alpha_1 \nu_N + \alpha_2\frac{\Delta x}{\Delta t} \nu_N^2+ \alpha_3 |\nu_N|\right) T^{-1}\\
		&= T\,\text{diag}(d(\nu_1), \ldots, d(\nu_N))\,T^{-1},
	\end{align*}
	
	Here, $\nu_i=\lambda_i\Delta t/\Delta x$ is the dimensionless characteristic speed and $d(\nu)$ is the dimensionless scalar dissipation function.
	
	Throughout the whole paper, $d(\nu)$ denotes the scalar, non-dimensional dissipation function. It can be understood as the effect of the dissipation matrix $D$ on the eigenvalues of $\tilde{A}$, as shown in Eq.~\eqref{eq:trafoDd}. 
	%
%==============================
\section{Classical Riemann Solvers}\label{sec:riemannsolvers}
%==============================
%
	In this section we will recall some well-known Riemann solvers, which are necessary for the development of the new family of hybrid Riemann solvers. The numerical flux function \eqref{eq:gerneralNumFlux} and thus the resulting scheme is completely dictated by the dissipation matrix $D$. As indicated by Eq.~\eqref{eq:trafoDd}, we can thus break down the discussion to the dimensionless scalar dissipation function $d(\nu)$. We only compare the dissipation functions and do not explicitly state the numerical flux functions. However, several well-known numerical fluxes are specified in appendix \ref{sec:Appendix}.
	
	 In order to obtain $L^2$ stability, $d(\nu)\geq \nu^2$ is required. If we want to restrict ourselves to monotone schemes, the dissipation function has to lie above the absolute value function \cite{GodlewskiRaviart2013}. Below this bound, the scheme is non-monotone, which means that undesired oscillations may occur at discontinuities. In general we want the dimensionless scalar dissipation function to fulfill $d(\nu)\geq |\nu|$.\newline 
	
	For a linear system, where we can write $\partial_x f(U)=A\,\partial_x U$ with some matrix $A(U) \in \mathbb{R}^{N \times N}$, which might be constant, the dissipation matrix of the upwind or Godunov scheme (both denoted by UP) reads
	\begin{align}
		D_\text{up}=|A|\quad\Leftrightarrow\quad d_\text{up}(\nu)=|\nu|.
	\end{align}
	For non-linear systems with general flux function $f$, where no matrix $A$ can be found such that $\partial_x f(U)=A\,\partial_x U$ holds true, the Godunov scheme can be adapted by using a so-called Roe Matrix $\tilde A(U_L, U_R)$, satisfying (amongst other properties) $\tilde A(U_L, U_R) (U_R-U_L) = F(U_R) - F(U_L)$, cf. \cite{Roe1981}. The dissipation matrix of Roe's scheme is given by
	\begin{align}
		D_\text{Roe}=|\tilde A|\quad\Leftrightarrow\quad d_\text{Roe}(\nu)=|\nu|,
	\end{align}
	The upwind scheme and the Roe solver are complete Riemann solvers. Since their dissipation function is exactly the absolute value function, and recalling that every dissipation function below the absolute value function is non-monotone, we can conclude that the upwind or Roe solver has the minimal dissipation while still being monotone. Thus it is optimal in the sense of minimal dissipation.
	
	Now follows a list of monotone and incomplete solvers with decreasing dissipation, all of them having more dissipation than the upwind scheme.
	
	The dissipation function of the monotone Lax-Friedrichs (LF) scheme is
	\begin{align}
		D_\text{LF}=\frac{\Delta x}{\Delta t}\,I\quad\Leftrightarrow\quad d_\text{LF}(\nu)=1.
	\end{align}	
	A solver which decreases the dissipation is the Rusanov scheme \cite{Rusanov1961}, also referred to as local Lax-Friedrichs (LLF) scheme. It takes into account the globally fastest eigenvalue of the system
	\begin{align}
		D_\text{LLF}=\max\left(\bar{\lambda}(U_L), \bar{\lambda}(U_R)\right)\,I\quad\Leftrightarrow\quad d_\text{LLF}(\nu)=\max(|\nu_{\min}|,|\nu_{\max}|),
	\end{align}	
	where $\nu_{\min}$ and $\nu_{\max}$ are the (non-dimensionalized) globally fastest left and right traveling wave speeds, i.e.
	\begin{subequations}
		\begin{align}
			\nu_{\min} &= \frac{\Delta t}{\Delta x}\min(\lambda_{\min}(U_L), \lambda_{\min}(U_R)),\\
			\nu_{\max} &= \frac{\Delta t}{\Delta x}\max(\lambda_{\max}(U_L), \lambda_{\max}(U_R)).				
		\end{align}	 
		\label{eq:nuMinNuMax}
	\end{subequations}
	
	Harten, Lax and van Leer \cite{HLL1983} further decreased the amount of dissipation added to the system by considering the fastest and slowest waves of the system
\begin{subequations}
	\begin{align}
		D_\text{HLL} &= \frac{|\lambda_L|\,\lambda_R - |\lambda_R|\,\lambda_L}{\lambda_R-\lambda_L}\,I + \frac{|\lambda_R|-|\lambda_L|}{\lambda_R-\lambda_L}\,\tilde{A}\\
		\Leftrightarrow d_\text{HLL}(\nu) &= \frac{|\nu_L|\,\nu_R - |\nu_R|\,\nu_L}{\nu_R-\nu_L} + \frac{|\nu_R|-|\nu_L|}{\nu_R-\nu_L}\,\nu.
	\end{align}
\end{subequations}
Here, $\lambda_L=\lambda_\text{min}(U_L)$, $\lambda_R=\lambda_\text{max}(U_R)$, and $\nu_{L/R}=\lambda_{L/R}\Delta t/\Delta x$. Note that there is a difference between $\nu_{L/R}$ and $\nu_{\min/\max}$, Eq.~\eqref{eq:nuMinNuMax}.

\begin{figure}[t]
		\centering %
		\includegraphics[width=0.5\textwidth]{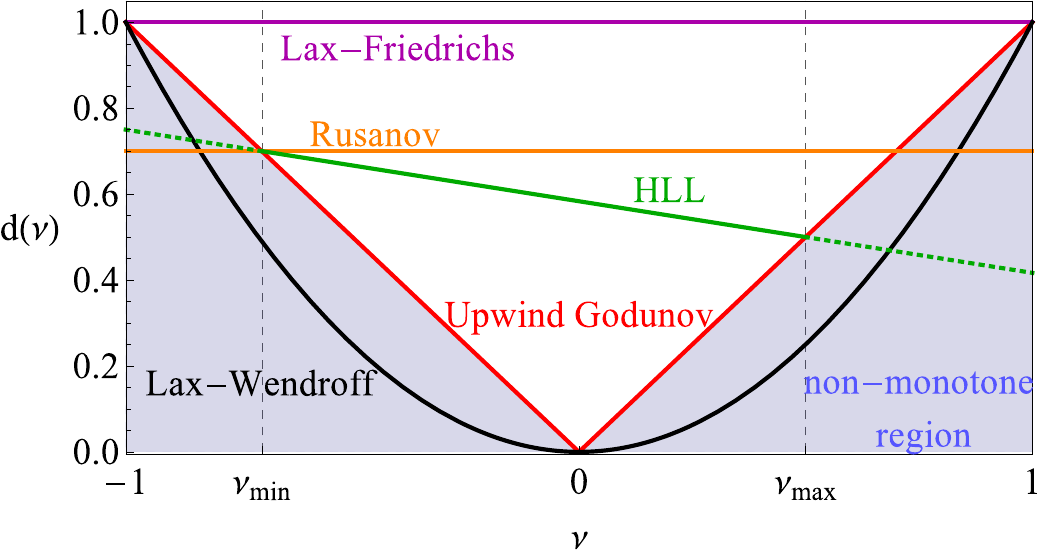}
		\caption{Comparing the scalar dissipation functions of different schemes.}
		\label{fig:differentSolvers}
	\end{figure}

A scheme which further reduces the dissipation is the Lax-Wendroff (LW) scheme, determined by
	\begin{align}
		D_\text{LW}=\frac{\Delta t}{\Delta x}\,\tilde A^2 \quad\Leftrightarrow\quad d_\text{LW}(\nu)=\nu^2.
	\end{align}	
	Discontinuities are approximated with steep gradients using the Lax-Wendroff scheme, however, the method is known to cause oscillations at discontinuities due to its non-monotonicity \cite{LeVeque1992}, see also Fig.~\ref{fig:differentSolvers}. 
%
%--------------------------------------
\subsection{Comparison of Classical Riemann Solvers}\label{subsec:compareRS}
%--------------------------------------
%
	Now that all classical solvers reviewed in Sec.~\ref{sec:riemannsolvers} are written in scalar, non-dimensional form, they can easily be compared. This is done in Fig.~\ref{fig:differentSolvers}, which shows the scalar dissipation of the methods (corresponding to the amount of dissipation added to the system) depending on $\nu \in \lbrack -1, 1\rbrack$. This restriction is due to the CFL condition (assuming an explicit Euler time integration) which requires that $\bar{\nu}=\bar{\lambda} \Delta t / \Delta x \leq 1$ and is fulfilled for all $\nu$ based on the eigenvalues of the flux Jacobian. During the whole paper we remain in the semi-discrete setting, thus, the choice of $\nu\in \lbrack -1, 1\rbrack$ is a conservative one, as other time integration schemes might allow for a wider range of $\nu$.
	
	 The shaded region denotes the non-monotone part. The upwind scheme marks the limit between monotone and non-monotone region. Thus, it is the scheme with the smallest amount of dissipation while still being monotone, i.e. not causing oscillations. The other extreme amongst the monotone schemes is the Lax-Friedrichs method, which induces more dissipation (because of larger values of $d(\nu)$) than all other methods considered in this paper. The Rusanov scheme reduces the dissipation. However, especially the slow waves, $\nu\approx 0$, suffer from too much dissipation. This is because for $\nu\approx 0$, the scalar dissipation function $d_{\text{LLF}}$ yields larger values than $d_{\text{up}}$ and therefore, more dissipation than needed is added to the system. This can also be seen in Fig.~\ref{fig:differentSolvers}. The Lax-Wendroff scheme on the other side of the spectrum has very little dissipation, however, it causes oscillations because it is non-monotone. Thus, concerning the balance between little dissipation and no spurious oscillations, the upwind scheme seems to be the scheme of choice. However, this scheme necessitates the decomposition of the eigensystem. In some cases no analytic form of the eigensystem is available but an estimate of the fastest and slowest eigenvalues can be obtained. Using only these two information, the HLL solver can be computed. Compared to solvers which only use the globally fastest wave speed, HLL reduces the amount of dissipation. Especially for faster waves close to the globally fastest signal velocities, HLL intersects with the upwind Godunov method, see Fig.~\ref{fig:differentSolvers}. Note that the Lax-Wendroff scheme is not the only scheme lying outside the monotone region. For the Rusanov and HLL schemes, the wave speeds have to be chosen  in such a way, that they bound the actual wave speeds.
%
%===========================================
%===========================================
\section{Hybrid Riemann Solvers}\label{sec:hybridSolvers}
%===========================================
%===========================================
%
	In this section, we introduce some hybrid Riemann solvers. This means, solvers which can be constructed using weighted combinations of the classical solvers in Sec.~\ref{sec:riemannsolvers}. %Our aim is to construct a solver, which requires as little information as HLL but has less dissipation. 
	We are seeking for solvers which require as little information as HLL but are less dissipative.  This is advantageous, especially because slow waves will be better resolved. Additionally, the solver should be computationally not more or only marginally more expensive than HLL. As an input, we only require the knowledge (or an estimate) of the globally slowest and fastest characteristic waves of the system, $\lambda_\text{min}$ and $\lambda_\text{max}$, and not the complete spectrum of eigenvalues.	
%
%--------------------------------------
\subsection{FORCE}\label{subsec:FORCE}
%--------------------------------------	
%
	We first consider the First Order Centered (FORCE) scheme, introduced by Toro et al. \cite{ToroBillett1996, ToroBillett2000}. This scheme can be viewed as a monotone version of the Lax-Wendroff method and the numerical flux function can be expressed as the average of the Lax-Friedrichs and the Lax-Wendroff method. Consequently, the dissipation matrix and function are given by
	\begin{align}
		D_\text{FORCE} = \frac{1}{2} \frac{\Delta x}{\Delta t} \left( \frac{\Delta t^2}{\Delta x^2}\,A^2 + I \right)\quad \Leftrightarrow\quad d_\text{FORCE}(\nu) =\frac{1}{2} \left( \nu^2 + 1 \right).
		\label{eq:FORCE}
	\end{align}	
	Note that this solver needs as little characteristic information as the Lax-Friedrichs scheme while needing one more flux evaluation.
%
%--------------------------------------
\subsection{MUSTA}\label{subsec:MUSTA}
%--------------------------------------	
%
In \cite{Toro2006} and \cite{ToroTitarev2006}, a multi stage flux called MUSTA has been been introduced. This flux function is based on the repetition of a simple flux function in order to resolve the Riemann fan and obtain high-resolution solutions of the Riemann problem. The general $\text{MUSTA}_k$ flux consists of $k$ repetitions. Here, we only consider the schemes for $k\in \{0, 1\}$. It turns out that $\text{MUSTA}_0$ is the FORCE scheme, Eq.~\eqref{eq:FORCE}, as suggested in \cite{Toro2006, ToroTitarev2006}. For $k=1$ the dissipation matrix of this scheme is given by
\begin{align}
	D_{\text{MUSTA}_1} = \frac{1}{4}\frac{\Delta x}{\Delta t}\,I +  \frac{\Delta t}{\Delta x}\,\tilde{A}^2 - \frac{1}{4}\left(\frac{\Delta t}{\Delta x}\right)^3\,\tilde{A}^4 \quad
	\Leftrightarrow \quad d_{\text{MUSTA}_1} = \frac{1}{4} + \nu^2 - \frac{1}{4}\nu^4.
\end{align}
It is interesting to realize that $d_{\text{MUSTA}_1}$ slightly drops below the absolute value function for larger wave speeds, see Fig.~\ref{fig:hybridSolvers}. This feature will be discussed in more detail in Sec.~\ref{subsec:Domega}.
%
%--------------------------------------------------------------------------------
\subsection{HLLX} \label{subsec:HLLX}
%--------------------------------------------------------------------------------
%
In this section we introduce a solver which also includes a quadratic term in the dissipation matrix. This means, it has the same number of flux evaluations as FORCE and less than $\text{MUSTA}_k$ for $k\geq 1$. 
%
%--------------------------------------------------------------------------------
\subsubsection{The $P_2$-Dissipation Function} \label{subsec:P2}
%--------------------------------------------------------------------------------
%
	%
	\begin{figure}[t]
		\centering %
		\includegraphics[width=0.5\textwidth]{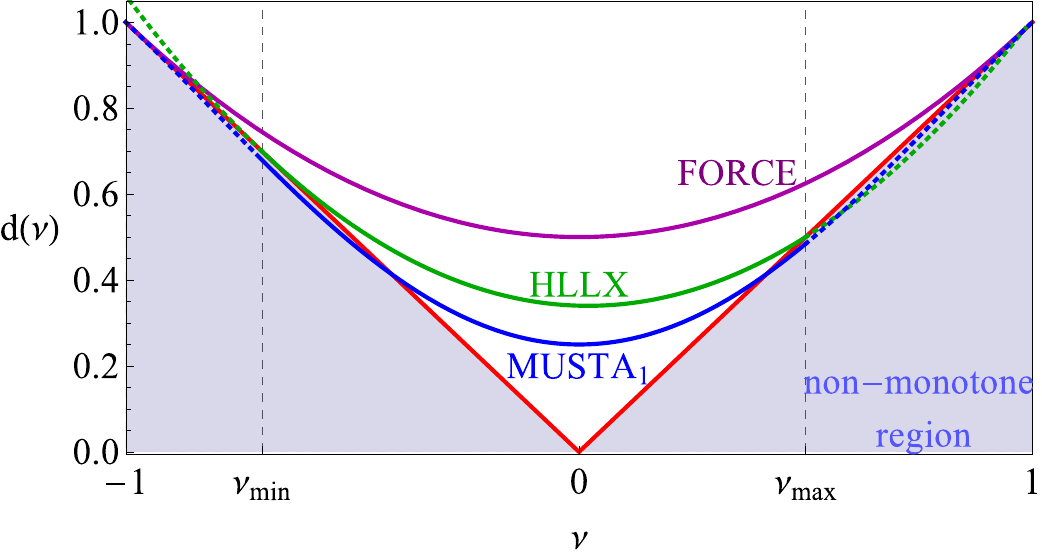}
		\caption{Scalar dissipation functions of hybrid Riemann solvers FORCE, $\text{MUSTA}_1$ and HLLX alias $P_2$.}
		\label{fig:hybridSolvers}
	\end{figure}
	 We recall three requirements which have been proposed by Degond et al. \cite{DegondPeyrardRussoVilledieu1999}. The resulting monotone Riemann solver, named $P_2$, is based on a quadratic dissipation function, fully determined by 
\begin{subequations}
	\label{eq:P2requirements}
	\begin{align}
		d_{P_2}(\nu_\text{min})&=d_\text{up}(\nu_\text{min})=|\nu_\text{min}|,\\
		d_{P_2}(\nu_\text{max})&=d_\text{up}(\nu_\text{max})=|\nu_\text{max}|,\\
		\text{and}\;\;d_{P_2}'(\bar\nu)&=d_\text{up}'(\bar\nu)=\text{sign}(\bar\nu), \quad\bar\nu=\begin{cases}\nu_\text{max}\;\,{if}\;\;|\nu_\text{max}|\geq|\nu_\text{min}|\\
			\nu_\text{min}\;\,\text{if}\;\;|\nu_\text{min}|>|\nu_\text{max}|.	\end{cases}
	\end{align}
\end{subequations}		 
	 The first two requirements indicate that the dissipation matrices match the absolute value function at the minimal and maximal Courant number. This means that the numerical flux function at these points equals the upwind flux. The third requirement is that the slope at the globally fastest wave speed has to match the one of the absolute value function. The dissipation function automatically fulfills $d_{P_2}(\nu)\geq |\nu|$ for $\nu \in [\nu_\text{min}, \nu_\text{max}]$, which means that it is monotone in this region. Even though the requirements \eqref{eq:P2requirements} have been proposed in \cite{DegondPeyrardRussoVilledieu1999}, neither the dissipation matrix nor the numerical flux function of $P_2$ have been explicitly stated.
%
%--------------------------------------------------------------------------------
\subsubsection{The HLLX-Dissipation Function} \label{subsec:DHLLX}
%--------------------------------------------------------------------------------
%	 
	We now present a simple way of implementing $P_2$, based on the flux functions of Lax-Friedrichs, HLL, and Lax-Wendroff, and call this Riemann solver HLLX. First, we notice that the dissipation function consists of a constant, an affine linear, and a quadratic part. Thus, it can be expressed as the weighted average 
	\begin{subequations}
		\begin{align}
		\label{eq:dissipFuncHLLX}
		d_\text{HLLX}(\nu) &=\alpha_0\,d_\text{LF}(\nu) +\alpha_1\,d_\text{HLL}(\nu)+\alpha_2\,d_\text{LW}(\nu),\\
		\intertext{where the coefficients are given by}
		\alpha_0 &= \alpha |\nu_{\min} \cdot\nu_{\max} |,\\
		\alpha_1 &= 1-\alpha(|\nu_{\max}|+|\nu_{\min}|), \\
		\alpha_2 &=\alpha, \\
		\text{and}\;\;\;\alpha &= \dfrac{\nu_{\max}-\nu_{\min}-\big| |\nu_{\max}|-|\nu_{\min}|\big|}{(\nu_{\max}-\nu_{\min})^2}. \label{eq:alpha}	
		\end{align}		
		\label{eq:dissipFuncHLLXwithCoeff}		
	\end{subequations}				
	An even simpler but not as demonstrative way of defining the dissipation function is
	\begin{align}		
		d_\text{HLLX}(\nu) &=d_\text{HLL}(\nu) + \alpha  (\nu-\nu_\text{min})(\nu-\nu_\text{max}),
		\label{eq:dP2alternative}
	\end{align}			
	with $\alpha$, as defined in Eq.~\eqref{eq:alpha}. 
	
	Based on Eq.~\eqref{eq:dP2alternative}, it is possible to write the dissipation matrix as
	\begin{align}
		D_\text{HLLX}(U_L, U_R) &=D_\text{HLL}(U_L, U_R) + \alpha \frac{\Delta t}{\Delta x} (\tilde{A}-\lambda_\text{min}I)(\tilde{A}-\lambda_\text{max}I).
	\end{align}
	This directly leads to the corresponding numerical flux function	
	\begin{align}
		\hat f_\text{HLLX}(U_L, U_R) &= \hat f_\text{HLL} -\frac{\alpha}{2} \frac{\Delta t}{\Delta x}(\tilde{A}-\lambda_\text{min}I)(\tilde{A}-\lambda_\text{max}I)(U_R - U_L).
	\end{align}
	%
%===========================================
%===========================================
\section{Jacobian-Free Implementation}\label{sec:JacobianFree}
%===========================================
%===========================================
%	
In the formulation of the numerical flux function $\hat{f}(U_L, U_R) = \frac{1}{2}\left(f(U_L) - f(U_R) \right) + \frac{1}{2} D(U_L, U_R)\; (U_R-U_L)$, it is necessary to know the flux Jacobian or Roe matrix $\tilde A$, in order to calculate the dissipation matrix $D(U_L, U_R)$. In case of the Euler equations, an explicit expression of the flux Jacobian matrix $A$ is available, however, for larger systems of conservation laws this might be difficult to compute or implement. In these cases, a Jacobian-free implementation is desirable. We now describe the implementation for the cases treated in this article, which are scaled versions of $\tilde A$ and $\tilde A^2$. For a generalization of the Jacobian-free implementation, the interested reader is referred to \cite[Sec.~3.2]{Torrilhon2012}.

In the flux formula \eqref{eq:gerneralNumFlux}, the dissipation matrix is always multiplied by the state difference $\Delta U = U_R-U_L$, so that the Jacobian matrix always appears as a matrix vector multiplication with $\Delta U$. We exploit this fact by using the finite difference formulation $DF(\bar{U})\;\Delta U = \lim_{\varepsilon\to 0}(f(\bar{U}+\varepsilon\,\Delta U)-f(\bar{U}))/\varepsilon$ with the average value $\bar{U}=0.5\,(U_L+U_R)$. In the following, assuming $\varepsilon\ll 1$, we use the implementation
\begin{subequations}
	\begin{align}
		A\,\Delta U &= \frac{f(\bar{U}+\varepsilon\,\Delta U)-f(\bar{U})}{\varepsilon}\;, \\
		A^2\,\Delta U &= \frac{f\left(
		\bar{U} + \varepsilon \frac{(f(\bar{U} + \varepsilon\Delta U)-f(\bar{U}))}{\varepsilon}		
		\right)-f\left(\bar{U}\right)}
		{\varepsilon}.
	\end{align}
\end{subequations}	
Another option is to note that the way of writing the dissipation function \eqref{eq:dissipFuncHLLXwithCoeff} directly leads to a simple form for the numerical flux function. In order to implement $\hat f_\text{HLLX}$ we do not actually need to compute the dissipation matrix $D_\text{HLLX}$ based on Eq.~\eqref{eq:dissipFuncHLLX}. We rather use the fact, that most users already implemented Lax-Friedrichs, HLL, and Lax-Wendroff (see appendix \ref{app:classicalFluxFcts} for full formulations of the numerical flux functions). Defining $\bar{f} = 0.5 (f(U_L) + f(U_R))$, we can rewrite the numerical flux functions \eqref{eq:gerneralNumFlux} of these three classical Riemann solvers to extract their dissipation matrices. Altogether this yields the numerical flux function
	\begin{align}
		\label{eq:fHLLXversion1}	
		\hat f_\text{HLLX}(U_L, U_R) &=\bar f + \alpha_0\,(\hat f_\text{LF}-\bar f) +\alpha_1\,(\hat f_\text{HLL}-\bar f)+\alpha_2\,(\hat f_\text{LW}-\bar f) \\
		&=\hat f_\text{HLL} + \alpha \,\left\lbrace |\nu_{\min} \cdot\nu_{\max}| (\hat f_\text{LF}-\bar f) -(|\nu_{\max}|+|\nu_{\min}|)\,(\hat f_\text{HLL}-\bar f) + (\hat f_\text{LW}-\bar f)\vphantom{\frac{.}{.}} \right\rbrace,
		\label{eq:fHLLXversion2}		
	\end{align}		
	with $\alpha$ as in Eq.~\eqref{eq:alpha}. This formulation does not include the Jacobian matrix $\tilde A$ and can therefore be used for any (large) system of hyperbolic conservation laws.
	
	The form of Eq.~\eqref{eq:fHLLXversion1} and \eqref{eq:fHLLXversion2} will also be used for the solvers presented in the next sections. The numerical flux functions can be found in appendix \ref{sec:Appendix}. 
%	
%==========================================================================================
\section{HLLX$\omega$ - a Family of Hybrid Riemann Solvers}\label{sec:HLLXomega}	
%==========================================================================================
%	
This section, which is the core of this paper, presents a new family of Riemann solvers, called HLLX$\omega$. These new solvers are based on a quadratic dissipation function in the form of HLLX. Their dissipation functions shall be closer to the absolute value function, i.e. the upwind scheme, for all waves $\lambda_i$ of the hyperbolic system, and thus for all $\nu_i$. This means, that the family of solvers yields even less dissipation than HLLX. The new schemes require as little information as HLL and HLLX, namely the globally slowest and fastest characteristic waves of the system, $\lambda_\text{min}$ and $\lambda_\text{max}$. Additionally, HLLX$\omega$ only needs the same number of flux evaluations as HLLX, namely two. Since we do not want to increase neither the number of input information, nor the number of flux evaluations, we lower the dissipation function by a certain amount. This amount is described by a parameter $\omega\in [0,1]$, which determines the monotonicity behavior of the solver. For $\omega=0$ we recover the monotone HLLX solver, and for $\omega=1$, the non-monotone Lax-Wendroff solver. All intermediate members of the HLLX$\omega$ family are slightly non-monotone for a certain range of waves. However, we show in this section that under some mild assumptions, the results do not show spurious oscillations.
%
%------------------------------------------------------------
\subsection{Beyond Monotonicity}\label{subsec:Domega}
\begin{figure}[t]
		\centering
		\includegraphics[width=0.5\textwidth]{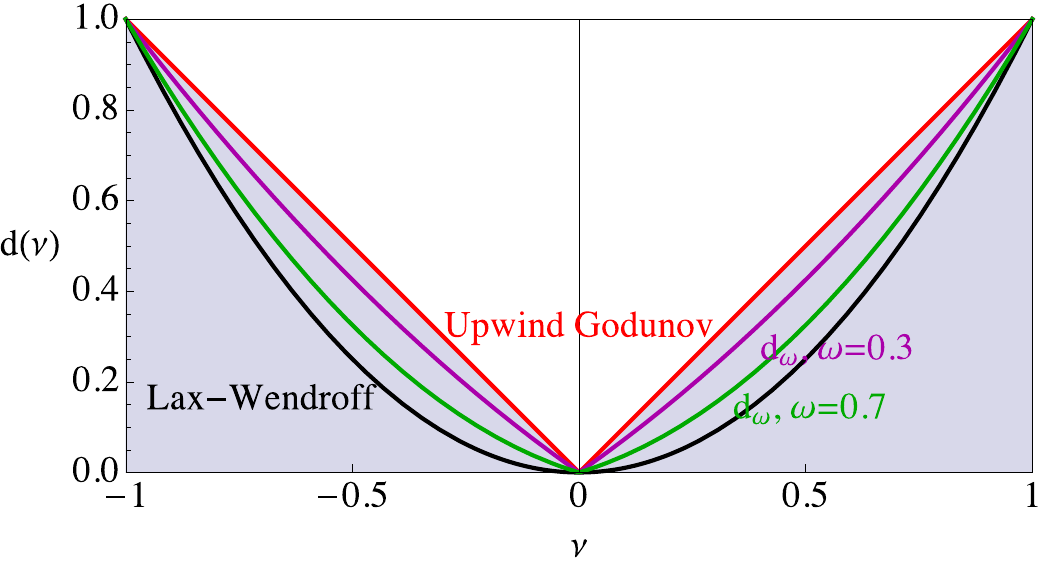}
		\caption{Scalar dissipation function $d_\omega(\nu)$, Eq.~\eqref{eq:dOmega} for $\omega=0.3$ and $\omega=0.7$.} 
		\label{fig:domega}
	\end{figure}

%------------------------------------------------------------
%	
Before we introduce this family of Riemann solvers, we state and validate some observations made by Torrilhon \cite{Torrilhon2012}. Firstly, it was perceived that the MUSTA fluxes introduced by Toro \cite{Toro2006} slightly drop below the upwind flux, which means that they do not fully lie in the monotonicity preserving region, see also Fig.~\ref{fig:hybridSolvers} and Sec.~\ref{subsec:MUSTA}. Therefore, as expected, the numerical solutions obtained with MUSTA fluxes show some non-monotone behavior. However, this behavior is far from the oscillations created by the Lax-Wendroff scheme. Additionally, the observed oscillations of MUSTA solutions decay in time and disappear after a certain number of time steps, cf. \cite[Fig. 5, p. A2084]{Torrilhon2012}. These results are essentially independent of the grid size and were observed for dissipation functions which only slightly drop below the absolute value function. Let us introduce the dissipation function $d_\omega(\nu)$, which is the weighted average of the dissipation functions of the monotone upwind scheme and the non-monotone Lax-Wendroff scheme,
	\begin{subequations}
		\begin{align}	
			d_\omega(\nu) &=\omega\ d_{\text{LW}}(\nu)+(1-\omega)\ d_{\text{UP}}(\nu) \\
			\Leftrightarrow\;\, D_\omega(U_L, U_R) &=\omega\ D_{\text{LW}}(U_L, U_R)+(1-\omega)\ D_{\text{UP}}(U_L, U_R),\;\; \omega\in [0,1].			
		\end{align}	
		\label{eq:dOmega}
	\end{subequations}	
	\hspace{-0.175cm} For $\omega=0$ we recover the monotone upwind scheme $d_{\omega=0}(\nu) = d_{UP}(\nu)$ and for $\omega=1, d_{\omega=1}(\nu) = d_{LW}(\nu)$ holds true. This can be seen in Fig.~\ref{fig:domega}, which shows $d_\omega(\nu)$ for $\omega\in \{0.0, 0.3, 0.7, 1.0\}$. 
	
	The aim of this section is to study the monotonicity behavior of $d_\omega(\nu)$ and produce similar effects to those found in \cite{Torrilhon2012}. We therefore investigate solutions of the scheme based on Eq.~\eqref{eq:dOmega} for different values of $\omega$. 
	We consider the scalar advection equation $u_t + u_x=0$ and use as initial condition the sign function, sgn$(x)$ on the interval $x\in [-1, 1]$. The jump evolves with time on a grid with $n=200$ grid cells until $T_\text{end}=0.25$. The CFL number is set to $\bar \nu=0.5$, which shows the maximal deviation of Lax-Wendroff from Upwind, cf. Fig~\ref{fig:domega}. For different values of $\omega$, we analyze the test case with the flux function resulting from Eq.~\eqref{eq:gerneralNumFlux} with \eqref{eq:dOmega}. 
	
\begin{figure}[t]
		\centering
		\begin{subfigure}{.4\textwidth}
			\includegraphics[width=\textwidth]{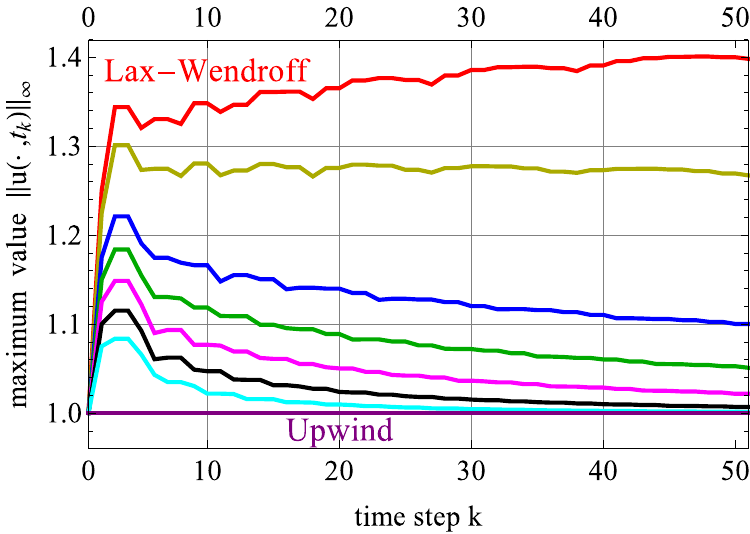}	
			\caption{Maximum value of the solution $u$ as a function of the number of time steps.}
			\label{fig:maxDomega}
		\end{subfigure}	
		\hfill
		\begin{subfigure}{.56\textwidth}
			\includegraphics[width=\textwidth]{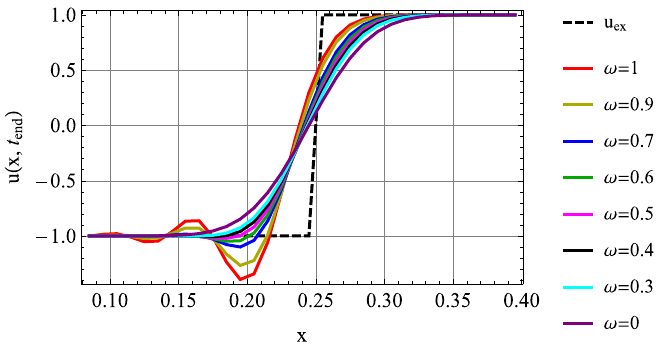}
			\caption{Zoom of the solution for different $\omega$.}
			\label{fig:uOmegas}
		\end{subfigure}	
		\caption{Test case with the sign function as initial condition on $x\in [-1, 1]$ with $n=200$ grid cells, $\,CFL= 0.5,$ and end-time $T_{\text{end}}=0.25$.}
		\label{fig:maxAndUdOmega}
	\end{figure}

	The numerical results for all tested values of $\omega$ are shown in Fig.~\ref{fig:uOmegas}. It can be easily seen, that $\omega=0$ and $\omega=1$ correspond to the upwind and the Lax-Wendroff schemes. That is, for $\omega=1$ we can observe the well-known oscillations. As $\omega$ decreases, the oscillations also decrease and for $\omega \lesssim 0.5$ no oscillations seem to be apparent in Fig.~\ref{fig:uOmegas}. For further analysis, the maximum value of the solutions $u$ as a function of the number of time steps is shown in Fig.~\ref{fig:maxDomega}, where 50 time steps correspond to $T_{\text{end}}=0.25$. Here, it can be seen that the oscillations of the Lax-Wendroff scheme do not decrease over time, whereas the Upwind scheme does not oscillate at any time. For the mixed schemes, a certain amount of weight needs to be given to the monotone Upwind scheme to make sure that the oscillations decrease in time. Now it can be seen what was not clearly visible in Fig.~\ref{fig:uOmegas}, namely that oscillations only completely decrease within 50 time steps for $\omega \to 0$. Note that these results are essentially independent of the grid size. 
	
	For the rest of the article we will assume that $\omega \leq 0.4$ is sufficiently small to diminish oscillations to an amount which can be considered "vanished". This assumption is tested and verified in the numerical experiments in Sec.~\ref{sec:numericalresults}.
	%
%----------------------------------------------------------------------------
\subsection{Modified Equation}\label{subsec:ModEq}
%----------------------------------------------------------------------------
%
	The phenomenon observed in Sec.~\ref{subsec:Domega} can be explained by taking a look at the modified equation of the scheme. This equation is obtained, when the difference equation of a numerical scheme is modeled by a differential equation \cite{LeVeque1992}. More specifically, the modified equation is the differential equation which is more accurately solved by the numerical scheme than the original equation \eqref{eq:conservation_law}. A scheme which solves Eq.~\eqref{eq:conservation_law} with order $p$, solves the modified equation
	\begin{align}
		\partial_t U(x,t)+\partial_x f(U(x,t)) &= D(\nu) \partial_{x}^{p+1} U(x, t)
		\label{eq:generalModEq}
	\end{align}
	with order $p+1$. Here, $D(\nu)$ is the dissipation matrix which can be obtained by computing the local truncation error of the method. Solving the linear advection equation $\partial_t U +a \partial_x U=0$ with the Upwind Godunov and the Lax-Wendroff scheme, their modified equations read
	\begin{subequations}
		\begin{align}
			\label{eq:modEqUP}
			\text{UP:}\quad\partial_t U(x, t)+a\,\partial_x U(x, t)  &= \frac{1}{2} a \Delta x (1-\nu) \partial_{xx} U(x, t)
			= D_\textsc{UP}(\nu)\, \partial_{xx} U(x, t), \\
			\text{LW:}\quad	\partial_t U(x, t)+a\,\partial_x U(x, t)  &= \frac{1}{6} a \Delta x^2 \left(\nu ^2-1\right) \partial_{xxx} U(x, t) 
			= D_\textsc{LW}(\nu)\, \partial_{xxx} U(x, t).
			\label{eq:modEqLW}
		\end{align}
		\label{eq:modEqUPandLW}
	\end{subequations}	
	Now we define the shift $\xi := x-a t$ and the shifted function $\tilde{u}(\xi, t)=:u(x, t)$, as well as the Fourier transform
	\begin{subequations}
		\begin{align}
			\label{eq:DefFourierTrafo}
			\hat{u}(k, t) = \mathcal{F}[\tilde{u}(\xi, t)] &:= \frac{1}{\sqrt{2 \pi }}\int _{-\infty }^{\infty } \tilde{u}(\xi, t) \exp(\mathbf{i} k \xi) d\xi,\\
			\tilde{u}(\xi, t) = \mathcal{F}^{-1}[\hat{u}(k, t)] &:= \frac{1}{\sqrt{2 \pi }}\int _{-\infty }^{\infty } \hat{u}(k, t)\exp(-\mathbf{i} k \xi) dk,
			\label{eq:DefInvFourierTrafo}
		\end{align}		
	\end{subequations}	
	using the default definition of modern physics. Applying the Fourier transform to Eq.~\eqref{eq:generalModEq} leads to an ordinary differential equation which can be solved in terms of $\hat{u}$ with solution
	\begin{align}
		\hat{u}(k, t) = \hat{u}_0(k) \exp(D(\nu)\,(-\mathbf{i} k)^{p+1}t).	
		\label{eq:uHatGeneral}
	\end{align}
	With the step function as initial condition, i.e. $\tilde{u}_0(\xi)=\sgn(\xi)$, we obtain $\hat{u}_0(k) = \sqrt{2/\pi}\;\mathbf{i}/k$, and thus
	\begin{align}
		\hat{u}(k, t) = \sqrt{\frac{2}{\pi }} \frac{1}{k}\, \mathbf{i}\, \exp(D(\nu)\,(-\mathbf{i} k)^{p+1}t).
		\label{eq:uHatWithIC}
	\end{align}
	The solution $\tilde{u}(\xi, t)$ can now be obtained by inverse Fourier transform \eqref{eq:DefInvFourierTrafo}, when inserting the diffusion coefficients $D_\textsc{UP}(\nu)$ and $D_\textsc{LW}(\nu)$. With the values described in Fig.~\ref{fig:maxAndUdOmega}, Eq.~\eqref{eq:modEqUPandLW} yields $D_\textsc{UP}=0.0025$ and $D_\textsc{LW}=-0.0000125$.
	In order to study the behavior of the scheme with diffusion matrix $D_\omega$ \eqref{eq:dOmega}, we compute its modified equation
	\begin{subequations}
		\begin{align}
			\partial_t U(x,t)+a\,\partial_x U(x,t)  &= \frac{1}{2} a \Delta x (1-\nu) (1-\omega)\partial_{xx} U(x, t) +\frac{1}{6} a\Delta x^2 \left(\nu ^2-1\right) \partial_{xxx} U(x, t)\\
			&= D_\textsc{UP}(\nu)(1-\omega)\, \partial_{xx} U(x, t) + D_\textsc{LW}(\nu)\, \partial_{xxx} U(x, t).
		\end{align}
		\label{eq:modEqDomega}
	\end{subequations}	
	\hspace{-0.175cm} For $\omega = 1$, the term $D_\textsc{UP}(\nu)(1-\omega)\, \partial_{xx} U(x, t)$ cancels and only the modified equation of Lax-Wendroff remains, as in Eq.~\eqref{eq:modEqLW}. For $\omega = 0$, we expect to recover the modified equation of the upwind scheme, Eq.~\eqref{eq:modEqUP}. This is the case including the second and third terms of the Taylor expansion of the local truncation error.

	\begin{figure}[t]
		\centering		
		\includegraphics[width=0.5\textwidth]{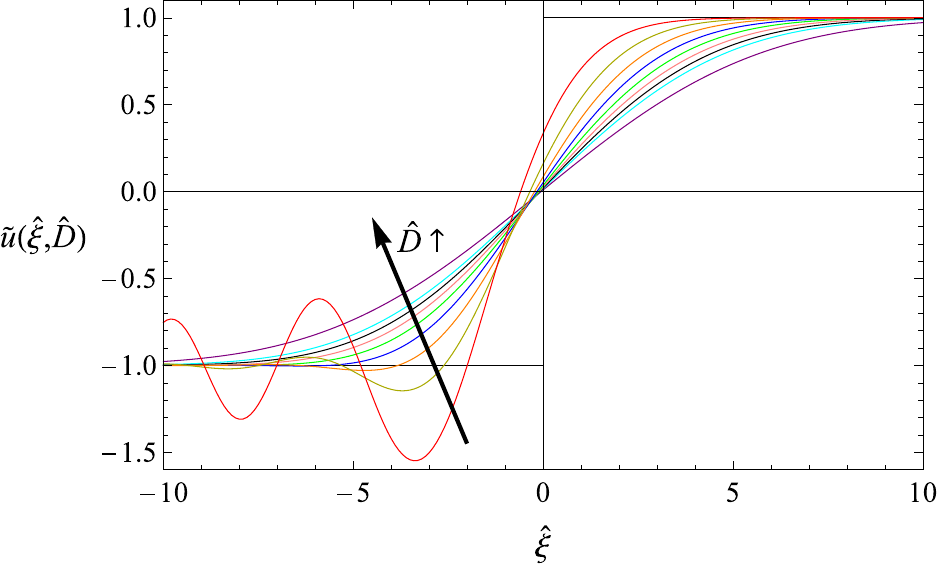}
		\caption{Analytical solution $\tilde{u}(\hat\xi, \hat{D})$ \eqref{eq:nonDimSoluhat}, depending only on the parameter $\hat{D}$ \eqref{eq:Dhat}.}
		\label{fig:analyticalSolutionWithDhat}
	\end{figure}

	Following the procedure described above, we obtain
	\begin{align}
		\hat{u}(k, t) = \sqrt{\frac{2}{\pi }} \frac{1}{k}\, \mathbf{i}\, \exp\left\lbrace(1-\omega)\;D_\textsc{UP} (-\mathbf{i}\,k)^2\;t + D_\textsc{LW} (-\mathbf{i}\,k)^3\;t\right\rbrace,\quad \omega\in[0, 1].
		\label{eq:uHatDomega}
	\end{align}
	Non-dimensionalizing Eq.~\eqref{eq:uHatDomega} leads to a formulation with only one remaining parameter. We set $\hat{k}=k/k_0$ and $\hat{\xi}=\xi/\xi_0$, with the constants
	\begin{align}
		k_0 = \sqrt[3]{\frac{1}{D_\textsc{LW}\,t}}, \quad\xi_0 = \frac{1}{k_0}=\sqrt[3]{D_\textsc{LW}\,t},\;\;\; \text{and }\;			\hat{D}= \sqrt[3]{\frac{t}{D_\textsc{LW}^2}}D_\textsc{UP}(1-\omega).
		\label{eq:Dhat}
	\end{align}
	Applying the inverse Fourier transform \eqref{eq:DefInvFourierTrafo} yields the solution
	\begin{align}
		\tilde{u}(\hat\xi, \hat{D}) = \frac{1}{\pi} \int _{-\infty }^{\infty } \frac{1}{\hat{k}} \exp(-\hat{D} \hat{k}^2)\left (\cos(\hat{k}^3) \sin(\hat{k} \hat{\xi}) - \sin(-\hat{k}^3) \cos(\hat{k} \hat{\xi})\right)\,d\hat{k},
		\label{eq:nonDimSoluhat}
	\end{align}
	which only depends on the non-dimensional parameters $\hat{\xi}$ and $\hat{D}$. Plotting $\hat\xi- \tilde{u}(\hat\xi, \hat{D})$ means that only the parameter $\hat{D}$ remains to yield different solutions, as depicted in Fig.~\ref{fig:analyticalSolutionWithDhat}. Increasing $\hat{D}$ can imply two things: 1) at a fixed value $\omega\in[0, 1]$, the time $t$ has been increased, corresponding to Fig.~\ref{fig:maxDomega}. Or 2), at a fixed time $t$, the parameter $\omega$ has been decreased, which corresponds to Fig.~\ref{fig:uOmegas}. Both scenarios result in smoother solutions and the disappearance of oscillations, as shown in Fig.~\ref{fig:analyticalSolutionWithDhat}.	
%
%-----------------------------------------------------
\subsection{HLL$\omega$}\label{subsec:HLLomega}
%-----------------------------------------------------
%
	Based on the findings of Sec.~\ref{subsec:Domega} and \ref{subsec:ModEq}, let us define a Riemann solver, called HLL$\omega$, which is a modification of HLL with less dissipation. Instead of intersecting with the absolute value function at $\nu_\text{min}$ and $\nu_\text{max}$, HLL$\omega$ fulfills the following constraints:
	\begin{align}
		d_{\text{HLL}\omega}(\nu_{\min})=d_\omega(\nu_{\min}), \quad d_{\text{HLL}\omega}(\nu_{\max})=d_\omega(\nu_{\max}).
	\end{align}
	These conditions yield the scalar dimensionless dissipation function 
	\begin{align}
		d_{\text{HLL}\omega}(\nu)=b_0+b_1 \nu.
	\end{align}
	The coefficients $b_0$ and $b_1$, as well as the numerical flux function $\hat f_{\text{HLL}\omega}(U_L, U_R)$ are given in appendix \ref{app:HLLomega}.
	
	The dissipation function $d_{\text{HLL}\omega}(\nu)$ is shown in Fig.~\ref{fig:HLLXomega}, where it is well-visible, that $\text{HLL}\omega$ yields solutions with less dissipation than HLL. Note that HLL$\omega$ is not monotone for all - however for most - wave speeds.
%
%-------------------------------------------------
\subsection{HLLX$\omega$}\label{subsec:HLLXomega}
\begin{figure}[t]
		\centering
		\includegraphics[width=0.5\textwidth]{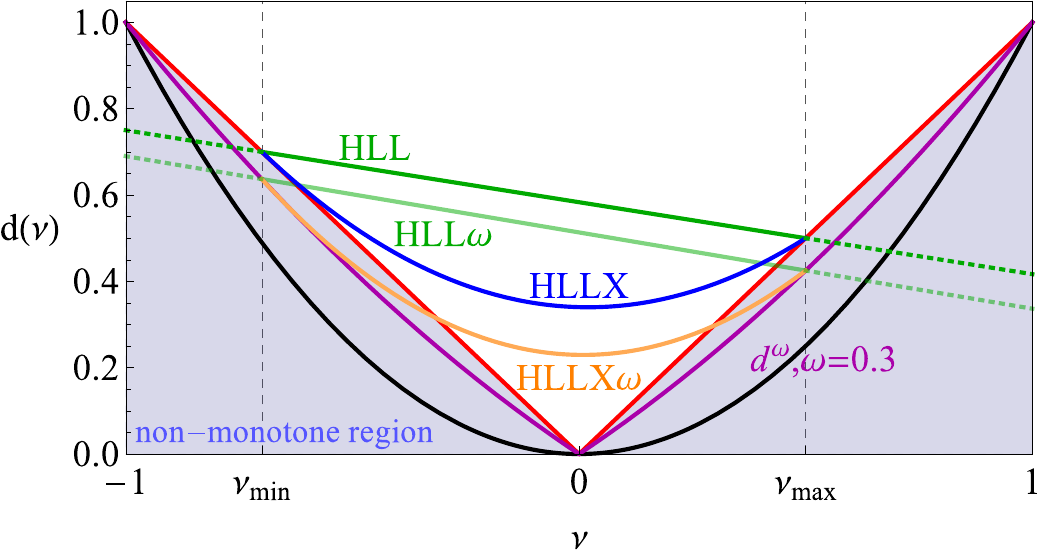}
		\caption{The monotone Riemann solvers HLL and HLLX and their non-monotone, however, less dissipative versions HLL$\omega$ and HLLX$\omega$. Here, $\omega=0.3$ has been chosen according to the findings of Sec.~\ref{subsec:Domega}.}
		\label{fig:HLLXomega}
	\end{figure}

%-------------------------------------------------
%
	Now we can come back to the aim of this section, the construction of a new family of approximate Riemann solvers, called HLLX$\omega$. Here, $\omega\in [0, 1]$ is a parameter, which controls the amount of dissipation of the solver. The improvement of HLLX$\omega$ is that the dissipation functions are closer to the absolute value function for all emerging wave speeds of the system and thus the solvers are less dissipative than HLLX. The dissipation functions of HLLX$\omega$ are designed in a similar fashion as $d_\text{HLLX}$, see Sec.~\ref{subsec:HLLX}. This means, $d_{\text{HLLX}\omega}(\nu)$ is a quadratic function, fully determined by 
	\begin{subequations}
		\begin{align}				 
			d_{\text{HLLX}\omega}(\nu_{\min})&=d_\omega(\nu_{\min}), \\
			d_{\text{HLLX}\omega}(\nu_{\max})&=d_\omega(\nu_{\max}), \\
			d_{\text{HLLX}\omega}'(\bar\nu)&=d_\omega'(\bar\nu),\;\; \text{where}\;\,\bar\nu=\begin{cases}\nu_\text{max},\quad\text{if}\;\;|\nu_\text{max}|\geq|\nu_\text{min}|\\
			\nu_\text{min},\quad\text{if}\;\;|\nu_\text{min}|>|\nu_\text{max}|.	\end{cases}
		\end{align}		
	\end{subequations}
	These conditions yield a dissipation function which can be written as a weighted linear combination of the Lax-Friedrichs, HLL$\omega$, and Lax-Wendroff dissipation functions or an extension of HLL$\omega$:
	\begin{subequations}
		\begin{align}					
			\label{eq:dHLLXomega}
			d_{\text{HLLX}\omega}(\nu) &= \beta_0\,d_\text{LF}(\nu) +\beta_1\,d_{\text{HLL}\omega}(\nu)+\beta_2\,d_\text{LW}(\nu),\\
			d_{\text{HLLX}\omega}(\nu) &= d_{\text{HLL}\omega}(\nu) + \beta (\nu-\nu_\text{min})(\nu-\nu_\text{max}).
			\label{eq:dHLLXomega2}			
		\end{align}		
	\end{subequations}
	The coefficients $\beta$ and $\beta_i$ can be found in appendix \ref{app:HLLXomega}. They depend on the parameter $\omega$, i.e. $\beta = \beta(\omega),\, \beta_i = \beta_i(\omega),\, i=0,1,2$. Therefore, the whole dissipation function of HLLX$\omega$ changes its behavior with $\omega$. Depending on the choice of this parameter, the dissipation function lies more or less inside the monotone region. For $\omega=0$, it is bound by the monotone HLLX solver, for $\omega=1$ HLLX$\omega$ recovers the $L^2$ stable but non-monotone Lax-Wendroff method.\\
	It can be seen in Fig.~\ref{fig:HLLXomega} that HLLX$\omega$ is less dissipative than HLL, HLL$\omega$, and HLLX. However, it does not fully lie in the monotone region, thus, one would expect oscillations near discontinuities. However, for the linear advection equation we observe that oscillations appearing close to discontinuities disappear after a certain number of time steps. For non-linear systems of equations, no oscillations are observed during the whole time of simulation. Thus, in any case, the final result obtained with HLLX$\omega$ is non-oscillatory, see discussion in Sec.~\ref{subsec:Domega} and \ref{subsec:ModEq}.

	As proposed for HLLX, there is an easy way of implementing the numerical flux function of HLLX$\omega$. By defining $\bar{f} = 0.5 (f(U_L) + f(U_R))$, we can rewrite the numerical flux function of Lax-Friedrichs, HLL$\omega$, and Lax-Wendroff to extract their dissipation matrices. These can then be used to formulate the Jacobian-free version of the numerical flux function
	\begin{align}
		\hat f_{\text{HLLX}\omega}(U_L, U_R) =\bar f + \beta_0\,(\hat f_\text{LF}-\bar f) +\beta_1\,(\hat f_{\text{HLL}\omega}-\bar f)+\beta_2\,(\hat f_\text{LW}-\bar f).
		\label{eq:HLLXomegaJacobianFree}
	\end{align}	
	Based on Eq.~\eqref{eq:dHLLXomega2}, it is also possible to write the dissipation matrix as
	\begin{align}
		D_{\text{HLLX}\omega}(U_L, U_R) &=D_{\text{HLL}\omega}(U_L, U_R) + \beta \frac{\Delta t}{\Delta x} (\tilde{A}-\lambda_\text{min}I)(\tilde{A}-\lambda_\text{max}I),\\
		\intertext{and thus, the numerical flux function can also be written as}
		\hat f_{\text{HLLX}\omega}(U_L, U_R) &= \hat f_{\text{HLL}\omega} +\frac{\beta}{2} \frac{\Delta t}{\Delta x}(\tilde{A}-\lambda_\text{min}I)(\tilde{A}-\lambda_\text{max}I)(U_L - U_R),
	\end{align}
	which is advantageous when the flux Jacobian $\tilde{A}$ is known.
	
	The choice of $\omega$ remains problem-dependent. However, $\omega\leq 0.5$ turned out to be a good choice and will be used in the following.
		
%==============================
\section{Numerical Results}\label{sec:numericalresults}
%==============================
%
In this section, we provide numerical experiments in order to demonstrate the performance of the new family of Riemann solvers described in Sec.~\ref{sec:HLLXomega}. As already stated, we are especially interested in large systems of conservation laws. Nevertheless, we start with the one-dimensional Euler equations with three emerging wave speeds. Already in this example, the difference in dissipation for slow waves (in this case the contact discontinuity) demonstrates different results. Other numerical examples are the ideal magnetohydrodynamics (MHD) equations, which exhibit seven characteristic velocities, and the 13-moment equations of Grad.

Since the focus of this paper is on the solvers themselves, all tests will be conducted with first order accuracy, using the explicit Euler method for time evolution.
% 
%--------------------------------------
\subsection{Sod's Shock Tube Problem}\label{subsec:Sod}
%--------------------------------------
%
Let us consider Sod's problem, which describes a shock tube containing two different ideal gases at the left and right side of a diaphragm, placed at $x=0$. The density, velocity, and pressure of the gases in the left and right region are given by
	\begin{align}
		\begin{pmatrix}
			\rho_L\\ \text{v}_L \\ p_L
		\end{pmatrix}
		=
		\begin{pmatrix}
			1.0\\ 0.0 \\ 1.0
		\end{pmatrix},\quad
		\begin{pmatrix}
			\rho_R\\ \text{v}_R \\ p_R
		\end{pmatrix}
		=
		\begin{pmatrix}
			0.125\\ 0.0 \\ 0.1
		\end{pmatrix}
	\end{align}
	At time $t=0$, the diaphragm is removed and the gases begin to mix. The time evolution is described by the one dimensional Euler equations, 
	\begin{figure}[t]
		\centering
		\begin{subfigure}{0.49\textwidth}
			\includegraphics[width=\textwidth]{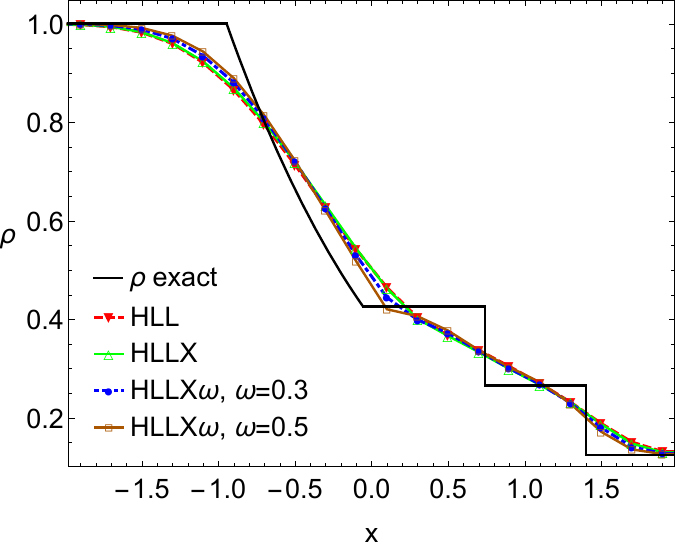}		
			\caption{Density profile.}		
			\label{fig:sodPbDensity}					
		\end{subfigure}		
		\begin{subfigure}{0.49\textwidth}
			\includegraphics[width=\textwidth]{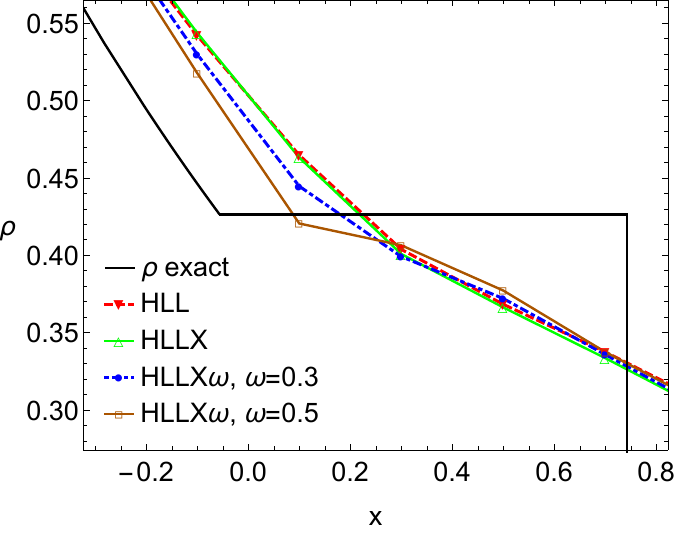}
			\caption{Zoom of the contact discontinuity.}
			\label{fig:sodPbZoomDensity}			
		\end{subfigure}		
		\caption{Solution of different Riemann solvers for Sod's shock tube problem on the domain $x\in[-2, 2]$, with  $N=20$ grid cells, $t_\text{end}=0.8, \bar\nu = 0.95$, and $\gamma=\frac{5}{3}$.}
		\label{fig:sodPb}
	\end{figure}	
	\begin{subequations}
		\begin{align}
			&\textbf{U}_t + \textbf{f(U)}_x = \textbf{0}
			\intertext{with the conserved variables \textbf{U}=$(\rho, \rho \text{v}, E)$, the flux function}
			&\textbf{f(U)}=\left(\rho \text{v}, \rho \text{v}^2 + p, \text{v}(E+p) \right)^T
			\intertext{and the equation of state for ideal gases}
			E=\frac{p}{\gamma-1}+\frac{1}{2}\rho \text{v}^2,
		\end{align}
		\label{eq:EulerEq}
	\end{subequations}
	\hspace{-.23cm} where the ratio of specific heats is set to $\gamma=1.4$. The computational domain is $[-2, 2]$ and the test is conducted with $N=20$ grid cells until $t_\text{end}=0.8$ with CFL number $\bar\nu = 0.95$.
	
	We compare the new family of solvers, HLLX$\omega$, for the two choices $\omega=0.3$ and $0.5$ with the HLL and HLLX solvers. The results for density are plotted in Fig.~\ref{fig:sodPb} together with the exact solution. Plots of velocity and pressure of the same simulation can be found in appendix \ref{app:Sod}. Sod's shock tube problem leads to three characteristic waves, where the pressure and velocity profile, Fig.~\ref{fig:sodPbPressure} and \ref{fig:sodPbVelocity}, illustrate the left-traveling rarefaction and the right-traveling shock wave. The contact discontinuity can only be seen in the density profile, Fig.~\ref{fig:sodPbDensity}. A zoom of the density profile around the contact discontinuity is depicted in Fig.~\ref{fig:sodPbZoomDensity}. This is of special interest since it highlights the different performances of the tested Riemann solvers because the contact discontinuity relates to the slow-moving wave of the system. As expected, HLLX$\omega$ with $\omega=0.5$ best approximates the steep gradient of the contact discontinuity since its dissipation function shows little diffusion around $\nu\approx 0$ compared to the other solvers. However, the solution shows slight overshoots in the velocity profile \ref{fig:sodPbVelocity}. HLLX$\omega$ with $\omega=0.3$ approximates the contact discontinuity almost as well as with $\omega=0.5$ and does not show spurious oscillations at $t_\text{end}=0.8$. HLLX and HLL are more dissipative than the new family of solvers, as expected from the discussions in Sec.~\ref{subsec:compareRS}, \ref{subsec:HLLX}, and \ref{sec:HLLXomega}.
%
%--------------------------------------
\subsection{Ideal Magnetohydrodynamics}
%--------------------------------------
%
	Ideal magnetohydrodynamics (MHD) describes the flow of plasma, assuming infinite electrical resistivity. The equations in one-dimensional processes read
	\begin{align}
	\label{eq:idealMHD}
		\partial_t \begin{pmatrix}
		\rho\\ \rho v_x\\ \rho\mathbf{v_t}\\ \mathbf{B_t}\\ E
	\end{pmatrix}	 
	+ \partial_x \begin{pmatrix}
		\rho v_x\\ \rho v_x^2+p+\tfrac{1}{2}\mathbf{B_t^2}\\ \rho v_x\mathbf{v_t}-B_x\mathbf{B_t}\\ v_x\mathbf{B_t}-B_x\mathbf{v_t}
		\\ (E+p+\tfrac{1}{2}\mathbf{B_t^2})v_x - B_x\mathbf{B_t}\cdot\mathbf{v_t}
	\end{pmatrix}
	 = 0
	\end{align}
	with density $\rho$, normal and transverse velocities $v_x$, and $\mathbf{v_t}=(v_y, v_z)$, respectively. Due to divergence constraints, the normal component of the magnetic field $B_x$ is constant in the one-dimensional case. The transverse magnetic field is $\mathbf{B_t}=(B_y, B_z)$. The energy $E$ is given in terms of the pressure $p$ by
	\begin{align}
		\label{eq:energy}
		E=\frac{1}{\gamma-1}p + \frac{1}{2}\rho (v_x^2 + \mathbf{v_t^2})+\frac{1}{2}\mathbf{B_t^2}.
	\end{align}
	The adiabatic constant $\gamma$ is set to $5/3$. System \eqref{eq:idealMHD} is hyperbolic and contains seven equations for the seven unknowns, $\mathbf{U}=(\rho, v_x,\mathbf{v_t}, p, \mathbf{B_t})$, exhibiting seven characteristic velocities. Thus, it can be considered as a large system of conservation laws.	

\begin{figure}[t]
		\centering
		\begin{subfigure}[t]{.48\textwidth}
			\includegraphics[width=\textwidth]{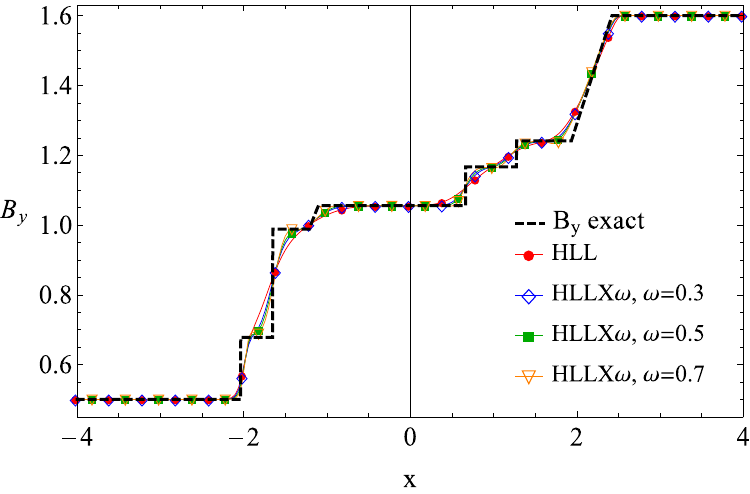}		
			\caption{Solution of the the magnetic field in $y$-direction.}					
		\end{subfigure}	
		\hfill
		\begin{subfigure}[t]{.49\textwidth}
			\includegraphics[width=\textwidth]{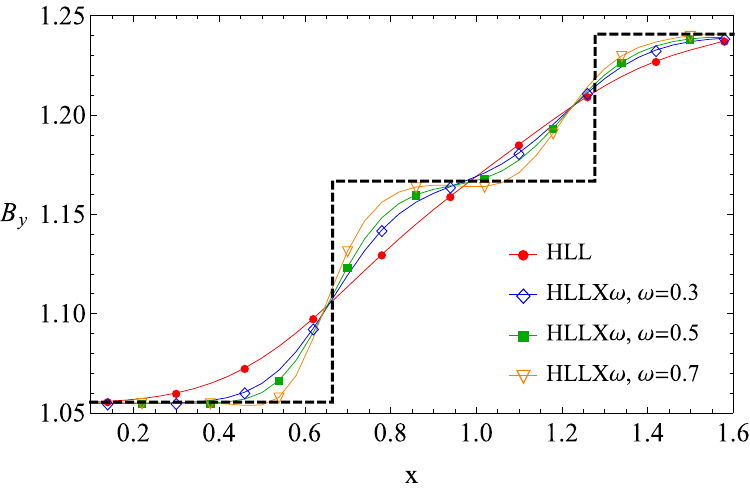}			
			\caption{Zoom of $B_y$ around the slow shock and rotational discontinuity to the right.}
			\label{fig:MHDsolutionZoomBy}
		\end{subfigure}	
		
		\begin{subfigure}[t]{.48\textwidth}
			\includegraphics[width=\textwidth]{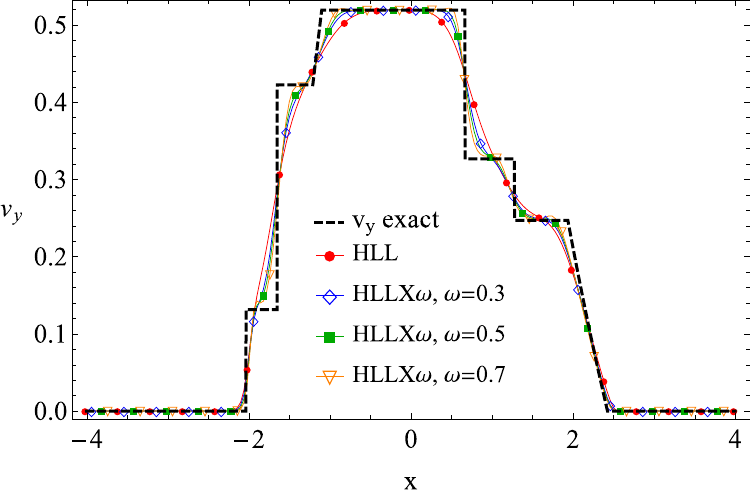}
			\caption{Solution of the the velocity in $y$-direction.}			
		\end{subfigure}	
		\hfill
		\begin{subfigure}[t]{.49\textwidth}
			\includegraphics[width=\textwidth]{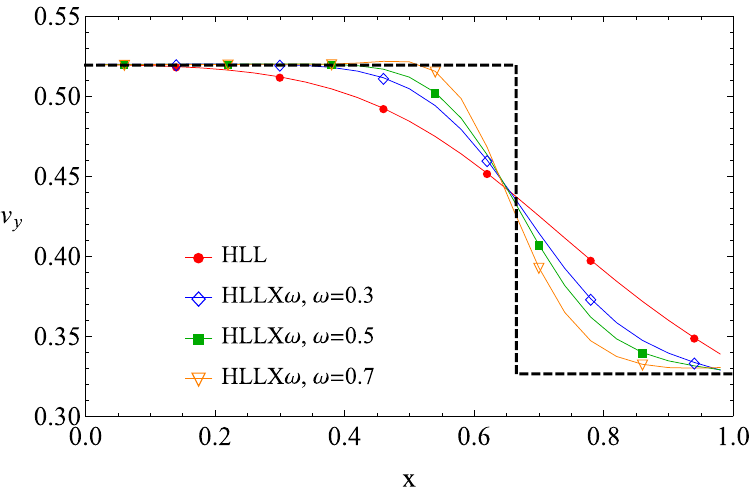}			
			\caption{Zoom of $v_y$ around the right-moving slow shock.}			
			\label{fig:MHDsolutionZoomVy}			
		\end{subfigure}	
		\caption{Solution of the ideal MHD equations with initial conditions \eqref{eq:ICmhd} at time $t_\text{end}=1.0$. Displayed are the $y$-components of magnetic field and velocity. The test has been computed with $x\in [-4,4],\,N=200$ grid cells, and CFL number $\bar\nu = 0.95$.}
		\label{fig:MHDsolution}
	\end{figure}

	Let us consider the Riemann problem given by
	\begin{subequations}
		\label{eq:ICmhd}
		\begin{align}
		(\rho^L, v_x^L, \mathbf{v_t}^L, p^L, \mathbf{B_t}^L)&=(1,0,(0,0),1,(0.5,0.6))\quad\text{if}\;\; x<0,\\
							 (\rho^R, v_x^R, \mathbf{v_t}^R, p^R,\mathbf{B_t}^R)&=(1,0,(0,0),1,(1.6, 0.2))\quad\text{if}\;\;x\geq 0,
		\end{align}
	\end{subequations}	
	and the normal magnetic field $B_x\equiv 1.5$. This problem, first introduced in \cite{TorrilhonExactRS}, represents a magnetic shock tube, since density, pressure, and velocity are constant in the whole domain and all fluid movements are generated only by the difference in the magnetic field. The solution of magnetic field and velocity in $y$-direction at time $t_{\text{end}}=1.0$ are shown in Fig.~\ref{fig:MHDsolution}. The test has been computed in the domain $[-4,4]$ with $N=200$ grid cells and CFL number $\bar\nu = 0.95$.
	We compare the new family of Riemann solvers with the parameter choices $\omega\in\{0.3, 0.5, 0.7\}$ to the HLL solver and the exact solution, which has been obtained by \cite{TorrilhonExactRS,Torrilhon2002ETH, Torrilhon2003Uniqueness,Torrilhon2003Convergence}. The magnetic field, as well as the velocity in $y$-direction exhibit all  waves of the system, except the contact discontinuity. The six waves are: a fast shock or rarefaction, a rotational discontinuity, and a slow shock or rarefaction, all of them to the left and to the right. The rotational discontinuity is also called Alfv{\'e}n wave. The shocks and discontinuities can be seen in Fig.~\ref{fig:MHDsolution}. Here, Fig.~\ref{fig:MHDsolutionZoomBy} shows a zoom of the right-moving slow shock and the right rotational discontinuity of the $y$-component of the magnetic field. This extract shows that the exact solution is better approximated with increasing value of $\omega$, which corresponds to decreasing dissipation. However, the solution of HLLX$\omega$ with $\omega=0.7$ yields some oscillations close to steep gradients at time $t_\text{end}=1.0$, so that it might preferable to use HLLX$\omega$ with $\omega=0.5$. In summary, all three simulations with HLLX$\omega$ yield closer approximations of the exact solution than HLL, which in comparison introduces more diffusion.
%
%----------------------------------------------
\subsubsection{Efficiency Study}
\begin{figure}
	\centering
	\begin{subfigure}[t]{.48\textwidth}
			\includegraphics[width=\textwidth]{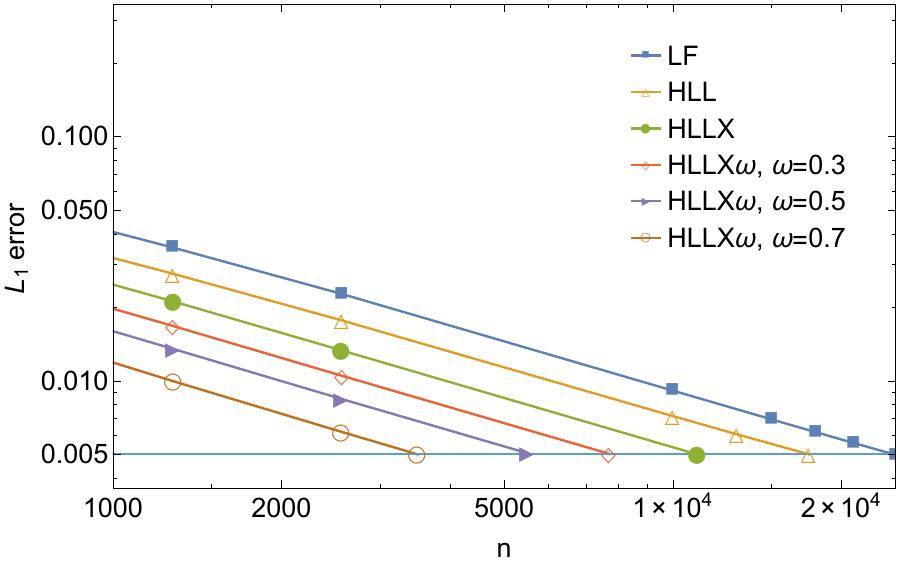}
			\caption{Grid refinement against $L_1$-error.}
			\label{fig:efficiencyNvsL1error}
		\end{subfigure}	
		\hfill
		\begin{subfigure}[t]{.49\textwidth}
			\includegraphics[width=\textwidth]{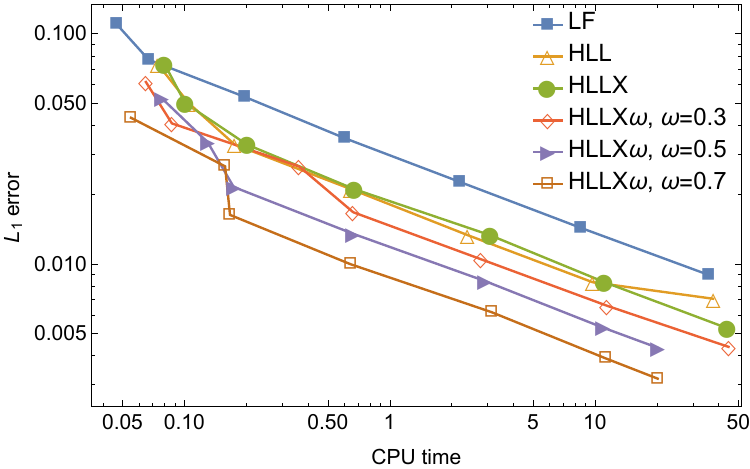}
			\caption{CPU time against $L_1$-error.}
			\label{fig:efficiencyCPUvsL1error}
		\end{subfigure}		
		\caption{Efficiency study of classical Riemann solvers and HLLX$\omega$.}	
\end{figure}

%----------------------------------------------
%	
	Finally, we illustrate that the new methods not only yield more accurate results but are also more efficient then classical methods. In Table \ref{tab:efficiency} we compare the performance of different schemes applied to the magnetic tube problem \eqref{eq:ICmhd}. For the five simulations, the time was measured until an $L_1$-error of $0.005$ was reached. The field "efficiency" stands for simulation time of the methods compared to the time of the HLL-simulation, i.e. $\textbf{t}_\text{HLL} / \textbf{t}_{\text{HLL}\mathbf{\star}}$. These results are also depicted in Fig.~\ref{fig:efficiencyNvsL1error}, which shows the mesh refinement and the corresponding error. For each method, the CPU time has been measured for the simulation corresponding to the mesh needed in order to obtain an $L_1$-error of $0.005$.
	
	Fig.~\ref{fig:efficiencyCPUvsL1error} shows the CPU time of different methods against the $L_1$-error corresponding to meshes with $N=20\cdot 2^j,\;j=3,\ldots,9$ grid points.
\begin{table}[h!]
\centering
	\begin{tabular}{ l | l }\hline
		Scheme & Efficiency \\ \hline
		HLLX$\omega$, $\omega =0.7$ & 22.396 \\
		HLLX$\omega$, $\omega =0.5$ & 8.057\\
		HLLX$\omega$, $\omega =0.3$ & 4.163\\
		HLLX &  2.253\\
		HLL &  1\\
		LF & 0.510\\ \hline
	\end{tabular}
	\caption{Performance table.}
	\label{tab:efficiency}
\end{table}
%
%--------------------------------------
\subsection{The 13-Moment Equations of Grad}\label{subsec:R13}
%--------------------------------------
%
	In order to show the performance of the new solver-family for a larger system of conservation laws, we consider the regularized 13-moment equations (R13). This fluid model describes rarefied fluids and micro-flows with high accuracy because it includes effects of higher moments. The R13 equations were derived from the Boltzmann equation by Struchtrup and Torrilhon \cite{StruchtrupTorrilhon2003} and are treated in detail in \cite{Torrilhon2006}. In this paper, we consider the homogeneous hyperbolic system of the R13 equations, i.e. source terms and gradient terms in the flux are neglected.
	
	The primitive variables of interest are mass density $\rho$, velocities $v_i$, pressure tensor $p_{ij}$, and heat flux $q_i$ with $i,j \in\{1,2,3\}$. The pressure tensor is symmetric, so that the R13 equations yield a total of 13 unknowns. 
	
	Using the Einstein notation, for the sake of simplicity, we define the total energy $E$ and the total energy flux $Q_i$
	\begin{subequations}		
		\begin{align}
	%		p &= \frac{1}{3}p_{ii}\quad \text{full pressure}\\
			E &= \frac{1}{2}\left( \rho v^2_k + p_{kk} \right),s\\
			Q_i &= q_i + E\,v_i + p_{ik}v_k
		\end{align}			
	\end{subequations}			
	and use the following tensor notation, $q_{(i} v_{j)} = \frac{1}{2}\left( q_i v_j + q_j v_i \right), \;		\delta_{(ij}q_{k)} = \frac{1}{3}\left( \delta_{ij}q_k + \delta_{jk}q_i + \delta_{ki}q_j \right)$.

	The one-dimensional, homogeneous R13 equations are given by
	\begin{align}
		\begin{cases}
			&\partial_t \rho + \partial_x(\rho v_1) = 0, \\
			&\partial_t \left(\rho v_i\right) + \partial_x (\rho v_1 v_i + p_{1i}) = 0,\quad i=1, 2, 3,\\
			&\partial_t \left(p_{ij}+\rho v_i v_j \right)+\partial_x\left(\rho v_1 v_i v_j + 3 p_{(ij} v_{1)} +\frac{6}{5}\delta_{(ij}q_{1)}\right) = 0,\quad i, j = 1, 2, 3,\\		
			&\partial_t Q_i + \partial_x\left( E v_i v_1 + 2 v_k p_{k(i} v_{j)} +\frac{2}{5}q_k v_k \delta_{i1}+\frac{14}{5}q_{(i} v_{1)} + \frac{1}{2}\left(p_{i1}v_k^2 + 5 \frac{p^2}{\rho}\delta_{i1}\right) \right) =0.
		\end{cases}	
	\end{align}		
	Let us consider the Riemann problem
	\begin{subequations}
		\label{eq:ICR13}
		\begin{align}
			\left(\rho^L, v_i^L, p_{ij}^L, q_i^L\right)_{i,j=1,2,3}
				&=\left(3, \begin{pmatrix}
				0 \\ 0.1\\ 0
				\end{pmatrix}, \begin{pmatrix}
				3 & 0 & 0 \\
				0 & 3 & 0 \\
				0 & 0 & 3
				\end{pmatrix}, \begin{pmatrix}
				0 \\ 0 \\ 0
\end{pmatrix}\right)\quad\text{if}\;\; x<0,\\
			\left(\rho^R, v_i^R, p_{ij}^R,q_i^R\right)_{i,j=1,2,3}
				&=\left(1, \begin{pmatrix}
				0 \\ 0\\ 0.1
				\end{pmatrix}, \begin{pmatrix}
				3 & 0 & 0 \\
				0 & 3 & 0 \\
				0 & 0 & 3
				\end{pmatrix}, \begin{pmatrix}
				0 \\ 0 \\ 0
\end{pmatrix}\right)\quad\text{if}\;\;x\geq 0.
		\end{align}
	\end{subequations}	
	The solution of density field and velocity in $y$-direction at time $t_{\text{end}}=0.9$ are shown in the appendix in Fig.~\ref{fig:R13density} and \ref{fig:R13velocity}. The test has been computed in the domain $[-2, 2.4]$ with $N=200$ grid cells and CFL number $\bar\nu = 0.8$.
	
	The new family of solvers with the parameter choices $\omega\in\{0.3, 0.5, 0.7\}$ is compared to the HLL and HLLX solvers. A reference solution has been computed on 3000 cells with the HLL flux. Again, the reference solution is better approximated with HLLX$\omega$ than with HLL or HLLX. Also, with increasing value of $\omega$, which corresponds to decreasing dissipation, the gradients of the numerical solutions are steeper. However, the solution of HLLX$\omega, \omega=0.7$ yields some oscillations, e.g. in density, Fig.~\ref{fig:R13density}, around $x\approx 1.1$, so that it might be preferable to use HLLX$\omega$ with $\omega=0.5$.
%	
%--------------------------------------
\subsection{Application to Higher-Order Methods in Space and Time}\label{subsec:higherOrder}
%--------------------------------------
%
The scope of this article is the development of the new family of Riemann solvers, HLLX$\omega$. Therefore, we performed all tests with first-order methods only. However, we did perform numerical simulations with second order reconstruction in space together with Heun's method. The obtained results were second order accurate which shows that the developed solvers can be used with higher order schemes in space and time. More test cases with other reconstruction techniques including limiters will be performed in the future.
%
%===========================================
\section{Conclusion}\label{sec:conclusion}
%===========================================

This paper presented a family of approximate hybrid Riemann solvers, HLLX$\omega$, for large non-linear hyperbolic systems of conservation laws.
The solvers do not require the characteristic decomposition of the flux Jacobian, only an estimate of the maximal propagation speeds in both directions is needed. The family of solvers contains a parameter $\omega$ which orders the solvers from fully-monotone to fully non-monotone. The intermediate solvers contain monotone as well as non-monotone parts. We showed that these intermediate family members, even though containing non-monotone parts for certain wave speeds, do not lead to oscillatory solutions in all test cases considered for this work. 

Extremely slow waves and stationary waves will still be approximated with higher dissipation than the upwind scheme, however, the computational cost of the new solvers is lower. Compared to solvers with similar prerequisites, the new Riemann solvers are able to rigorously decrease the dissipation of the scheme.

The numerical examples underline the excellent performance of the new family of solvers with respect to other solvers.
%
%==============================
\begin{appendix}
%==============================
\section{Appendix - Details on Riemann Solvers}\label{sec:Appendix}%Coefficients of Hybrid Riemann Solvers
%==============================
%
In this part of the appendix we state detailed formulations of numerical flux functions and dissipation functions/matrices of classical as well as the new Riemann solvers.
%
%------------------------------------------------------------
\subsection{Numerical Flux Functions of Classical Riemann Solvers}\label{app:classicalFluxFcts}
%------------------------------------------------------------
%
The numerical flux functions of the Lax-Friedrichs (LF), upwind (UP), HLL and Lax-Wendroff (LW) Richtmeyr schemes are given by \cite{LeVeque1992}
\begin{align}
	\hat f_{\text{LF}} &= \frac12\left(f(U_L)+f(U_R)\right) - \frac12\frac{\Delta x}{\Delta t}(U_R-U_L),\\
	\hat f_{\text{UP}} &= \frac12\left(f(U_L)+f(U_R)\right) - \frac{1}{2}|A|(U_R-U_L),
\end{align}	
\begin{align}
	\hat f_{\text{HLL}} &= \frac{\lambda_R f(U_L) − \lambda_L f(U_R) + \lambda_L\lambda_R(U_R − U_L)}{\lambda_R − \lambda_L} = \frac12\left(f(U_L)+f(U_R)\right) - \frac12 (a_0\,(U_R - U_L) + a_1\,(f(U_R)-f(U_L)) ),\\
	\hat f_{\text{LW}} &= f\left(\frac{1}{2}(U_L+U_R)-\frac{1}{2}\frac{\Delta t}{\Delta x}(f(U_R)-f(U_L))\right),
\end{align}
with the physical flux function $f$ specified by the hyperbolic conservation law and the coefficients \cite{Torrilhon2012}
\begin{align}
	a_0 = \frac{|\lambda_L|\,\lambda_R - |\lambda_R|\,\lambda_L}{\lambda_R-\lambda_L},\quad a_1 = \frac{|\lambda_R|-|\lambda_L|}{\lambda_R-\lambda_L}, \quad\lambda_L=\lambda_\text{min}(U_L),\quad \lambda_R=\lambda_\text{max}(U_R).
\end{align}
%
%------------------------------------------------------------
\subsection{HLL$\omega$}\label{app:HLLomega}
%------------------------------------------------------------
%
	As described in Sec.~\ref{subsec:HLLomega}, the scalar dissipation function of HLL$\omega$ is of the form $d_{\text{HLL}\omega}(\nu)=b_0+b_1 \nu$, where the coefficients $b_0 = b_0(\omega)$ and $b_1=b_1(\omega)$ are given by
	\begin{align}
		b_0(\omega) &= \dfrac{\nu_{\max}( \omega\,\nu_{\min}^2+(1-\omega) |\nu_{\min}|)-\nu_{\min}( \omega\, \nu_{\max}^2+(1-\omega) |\nu_{\max}|)}{\nu_{\max}-\nu_{\min}}, \\
		b_1(\omega) &= \dfrac{(1-\omega)(|\nu_{\max}|-|\nu_{\min}|)+\omega(\nu_{\max}^2-\nu_{\min}^2)}{\nu_{\max}-\nu_{\min}}.
	\end{align}
	Note that these coefficients are non-dimensional, which means that the dissipation matrix $D_{\text{HLL}\omega}$ (which has dimension $\frac{\Delta x}{\Delta t}$) reads
	\begin{align}
		D_{\text{HLL}\omega} = b_0(\omega)\,\frac{\Delta x}{\Delta t}\,I + b_1(\omega)\,A
	\end{align}
	and also $\lim_{\omega\to 0}b_0(\omega) = a_0\frac{\Delta t}{\Delta x}$ because $a_0$ has dimension $\frac{\Delta x}{\Delta t}$. In summary, the numerical flux function for HLL$\omega$ can be written as
	\begin{align}
		\hat f_{\text{HLL}\omega}(U_L, U_R) = \bar f - \frac12 \left( b_0(\omega) \frac{\Delta x}{\Delta t} (U_R - U_L) + b_1(\omega)(f(U_R) - f(U_L) \right)
	\end{align}	
	with $\bar{f} = 0.5 (f(U_L) + f(U_R))$.
%	
%------------------------------------------------------------
\subsection{HLLX$\omega$}\label{app:HLLXomega}
%------------------------------------------------------------
%
	The scalar dissipation function of HLLX$\omega$ can be written in different forms, see Sec.~\ref{subsec:HLLXomega}. The form with only one parameter is
	\begin{subequations}
		\begin{align}					
			d_{\text{HLLX}\omega}(\nu) &= d_{\text{HLL}\omega}(\nu) + \beta(\omega)\,(\nu-\nu_\text{min})(\nu-\nu_\text{max}), \vspace{-0.3cm}
			\intertext{where \vspace{-0.2cm}}
			\beta(\omega) &= \omega + (1-\omega)\frac{\nu_\text{max}-\nu_\text{min}-\big| |\nu_\text{max}| - |\nu_\text{min}| \big|}{(\nu_\text{max}-\nu_\text{min})^2}.
			\label{eq:betaCoeff}			
		\end{align}	
	\end{subequations}
	Note that $\beta(\omega) = \omega + (1-\omega)\alpha$, with the HLLX coefficient $\alpha$, Eq.~\eqref{eq:alpha}. Thus, it is easy to verify that for $\omega=0$, HLLX$\omega$ recovers the monotone HLLX solver.
	\newline\newline
	The dissipation function can also be written as a linear combination of Lax-Friedrichs, HLL$\omega$, and Lax-Wendroff
	\begin{subequations}
		\begin{align}					
			d_{\text{HLLX}\omega}(\nu) &= \beta_0(\omega)\,d_\text{LF}(\nu) +\beta_1(\omega)\,d_{\text{HLL}\omega}(\nu)+\beta_2(\omega)\,d_\text{LW}(\nu),
			\intertext{with the coefficients \vspace{-0.2cm}}
			\beta_0(\omega) &= \beta(\omega)\frac{(1-\omega)\left|\nu_{\min}\,\nu_{\max}\right|}{(1-\omega) + \omega  (\left|\nu_{\min}\right|+\left|\nu_{\max}\right|)}\;,\\
			\beta_1(\omega) &=1 - \beta(\omega)\left(\frac{1-\omega}{\left|\nu_{\min}\right|+\left|\nu_{\max}\right|}+\omega \right)^{-1},\\
			\beta_2(\omega) &=\beta(\omega),
		\end{align}	
	\end{subequations}
	with $\beta(\omega)$ given by Eq.~\eqref{eq:betaCoeff}. The coefficients $\beta$ and $\beta_i,\, i=0,1,2$ depend on $\omega \in [0,1]$. It is easy to verify that for $\omega=0$ the coefficients of HLLX are recovered, i.e.
	\begin{align*}
		\beta_0(\omega=0)=\alpha_0,\\
		\beta_1(\omega=0)=\alpha_1,\\
		\beta_2(\omega=0)=\alpha_2.
	\end{align*}
	This can also be seen in Eq.~\eqref{eq:dissipFuncHLLXwithCoeff} which describes the coefficients $\alpha_i$.
	
	The numerical flux function of HLLX$\omega$ can be written in the Jacobian-free formulation given by Eq.~\eqref{eq:HLLXomegaJacobianFree}
	\begin{align*}
		\hat f_{\text{HLLX}\omega}(U_L, U_R) =\bar f + \beta_0\,(\hat f_\text{LF}-\bar f) +\beta_1\,(\hat f_{\text{HLL}\omega}-\bar f)+\beta_2\,(\hat f_\text{LW}-\bar f).		
	\end{align*}
%==============================
\section{Appendix - Numerical Results}\label{sec:AppendixB}%Coefficients of Hybrid Riemann Solvers
%==============================
%
This part of the appendix yields more plots of solutions of the numerical test cases described in Sec.~\ref{sec:numericalresults}. This permits larger figures and therefore the performance of the different schemes becomes more clear.
%
%s------------------------------------------------------------
\subsection{Sod's Shock Tube Problem}\label{app:Sod}
%------------------------------------------------------------
%
Here we give the plots of velocity and pressure for Sod's shock tube problem as described in Sec.~\ref{subsec:Sod}. This means the solution has been computed on the domain $x\in[-2, 2]$, with  $N=20$ grid cells until time $t_\text{end}=0.8$. The CFL number is $\bar\nu = 0.95$, and the adiabatic coefficient $\gamma=\frac{5}{3}$.
	\begin{figure}[h!]
		\centering		
		\begin{subfigure}{0.49\textwidth}
			\includegraphics[width=\textwidth]{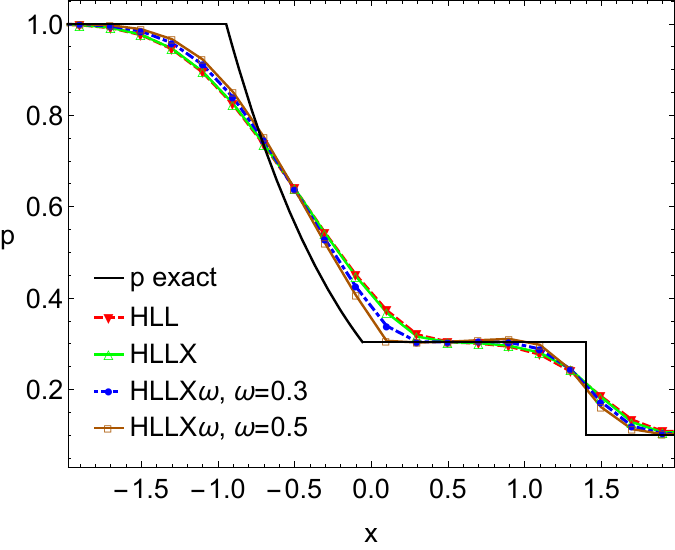}				
			\caption{Pressure profile.}
			\label{fig:sodPbPressure}
		\end{subfigure}		
		\begin{subfigure}{0.49\textwidth}
			\includegraphics[width=\textwidth]{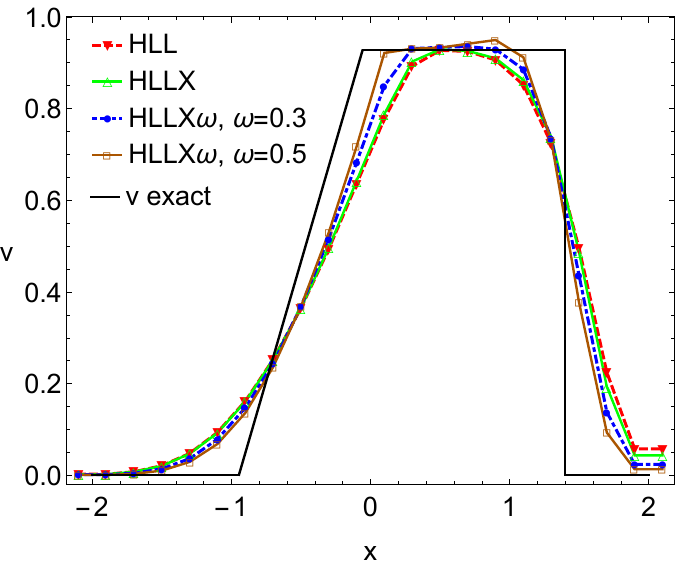}
			\caption{Velocity profile.}
			\label{fig:sodPbVelocity}			
		\end{subfigure}		
		\caption{Solution of different Riemann solvers for Sod's shock tube problem on the domain $x\in[-2, 2]$, with  $N=20$ grid cells, $t_\text{end}=0.8, \bar\nu = 0.95$, and $\gamma=\frac{5}{3}$.}
		\label{fig:sodPbApp}
	\end{figure}	
IIt can be seen that HLLX$\omega$ with $\omega=0.5$ best approximates the gradients of the discontinuities since the dissipation function has very little diffusion compared to the other solvers. However, the solution shows slight overshoots in the velocity profile \ref{fig:sodPbVelocity}. As already mentioned in Sec.~\ref{subsec:Sod}, HLLX$\omega$ with $\omega=0.3$ approximates the contact discontinuity almost as well as with $\omega=0.5$ and does not show spurious oscillations.
%
%------------------------------------------------------------
\subsection{The 13-Moment Equations of Grad}\label{app:R13}
\begin{figure}[h!]
		\centering
		\includegraphics[width=0.8\textwidth]{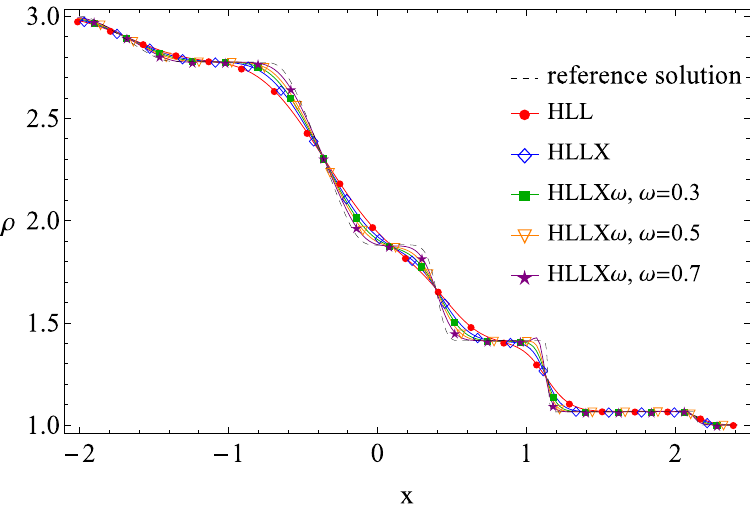}
		\caption{Solution of the of the R13 equations - density field.}
		\label{fig:R13density}
	\end{figure}

\begin{figure}[h!]
		\centering	
		\includegraphics[width=0.8\textwidth]{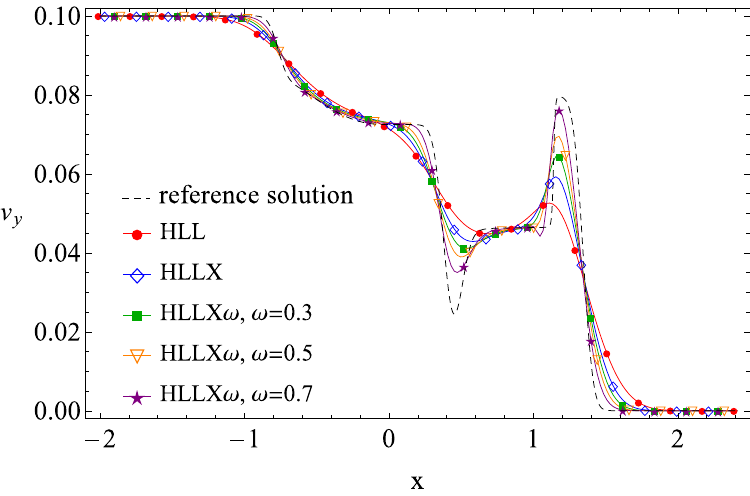}	
		\caption{Solution of the of the R13 equations - velocity field in $y$-direction.}
		\caption{}
		\label{fig:R13velocity}
	\end{figure}

%------------------------------------------------------------
%
Here we show the solution of the R13 equations described in Sec.~\ref{subsec:R13} with initial conditions \eqref{eq:ICR13}. The test has been computed with $x\in [-2,2.4],\,N=200$ grid cells and CFL number $\bar{\nu}= 0.8,\, T_{\text{end}}=0.9$. The reference solution has been obtained on 3000 cells with an HLL flux. Fig.~\ref{fig:R13density} shows the density field and \ref{fig:R13velocity} the velocity field in $y$-direction.
	%

%
%==============================
\end{appendix}
%==============================

%%%%%%%%%  REFERENCES %%%%%%%%
%=============================================================================================
\bibliographystyle{plain}
	\clearpage
	\bibliography{./references_complete}

\begin{thebibliography}{10}

\bibitem{CadaTorrilhon2009}
M.~Cada and M.~Torrilhon.
\newblock Compact third order limiter functions for finite volume methods.
\newblock {\em J. Comput. Phys.}, 228(11):4118--4145, 2009.

\bibitem{CastroGallardoMarquina2014}
M.~J. Castro, J.~M. Gallardo, and A.~Marquina.
\newblock A class of incomplete riemann solvers based on uniform rational
  approximations to the absolute value function.
\newblock {\em J. Sci. Comput.}, 60(2):363--389, 2014.

\bibitem{ColellaWoodward1984}
P.~Colella and P.~R. Woodward.
\newblock The piecewise parabolic method ({PPM}) for gas-dynamical simulations.
\newblock {\em J. Comput. Phys.}, 54:174--201, 1984.

\bibitem{DegondPeyrardRussoVilledieu1999}
P.~Degond, P.-F. Peyrard, G.~Russo, and P.~Villedieu.
\newblock Polynomial upwind schemes for hyperbolic systems.
\newblock {\em Comptes Rendus de l'Acad{\'e}mie des Sciences -- Series I --
  Mathematics}, 328(6):479--483, 1999.

\bibitem{GodlewskiRaviart2013}
E.~Godlewski and P.-A. Raviart.
\newblock {\em Numerical approximation of hyperbolic systems of conservation
  laws}, volume 118.
\newblock Springer Science \& Business Media, 2013.

\bibitem{HLL1983}
A.~Harten, P.~D. Lax, and B.~van Leer.
\newblock On upstream differencing and godunov-type schemes for hyperbolic
  conservation laws.
\newblock {\em SIAM review}, 25(1):35--61, 1983.

\bibitem{JiangShu1996}
G.-S. Jiang and C.-W. Shu.
\newblock Efficient implementation of weighted {ENO} schemes.
\newblock {\em J. Comput. Phys.}, 126(1):202--228, 1996.

\bibitem{LeVeque1992}
R.~J. Le~Veque.
\newblock {\em Numerical methods for conservation laws}.
\newblock Birkh{\"a}user, second edition, 1992.

\bibitem{LeVeque2002}
R.~J. Le~Veque.
\newblock {\em Finite volume methods for hyperbolic problems}.
\newblock Cambridge University Press, first edition, 2002.

\bibitem{LiuOsherChan1994}
X.-D. Liu, S.~Osher, and T.~Chan.
\newblock Weighted essentially non-oscillatory schemes.
\newblock {\em J. Comput. Phys.}, 115:200--212, 1994.

\bibitem{Marquina1994}
A.~Marquina.
\newblock Local piecewise hyperbolic reconstruction of numerical fluxes for
  nonlinear scalar conservation laws.
\newblock {\em SIAM J. Sci. Comput.}, 15:892--915, 1994.

\bibitem{Roe1981}
P.~L. Roe.
\newblock Approximate {R}iemann solvers, parameter vectors and difference
  schemes.
\newblock {\em J. Comput. Phys.}, 43:357--372, 1981.

\bibitem{Rusanov1961}
V.~V. Rusanov.
\newblock Calculation of interaction of non-steady shock waves with obstacles.
\newblock {\em J. Comput. Math. Phys. USSR}, 1:267–--279, 1961.

\bibitem{SchmidtmannSeiboldTorrilhon2015}
B.~Schmidtmann, B.~Seibold, and M.~Torrilhon.
\newblock Relations between {WENO3} and third-order limiting in finite volume
  methods.
\newblock {\em J. Sci. Comput.}, pages 1--29, 2015.

\bibitem{StruchtrupTorrilhon2003}
H.~Struchtrup and M.~Torrilhon.
\newblock Regularization of {G}rad's 13 moment equations: {D}erivation and
  linear analysis.
\newblock {\em Phys. Fluids}, 15:2668--2680, 2003.

\bibitem{ToroRiemannSolvers}
E.~F. Toro.
\newblock {\em Riemann solvers and numerical methods for fluid dynamics: a
  practical introduction}.
\newblock Springer, 2009.

\bibitem{Toro2006}
E.F. Toro.
\newblock {MUSTA}: A multi-stage numerical flux.
\newblock {\em Appl. Numer. Math.}, 56(10):1464--1479, 2006.

\bibitem{ToroBillett1996}
E.F. Toro and S.J. Billett.
\newblock Centred {TVD} schemes for hyperbolic conservation laws.
\newblock Technical report, Technical Report MMU–9603, Department of
  Mathematics and Physics, Manchester Metropolitan University, UK., 1996.

\bibitem{ToroBillett2000}
E.F. Toro and S.J. Billett.
\newblock Centred {TVD} schemes for hyperbolic conservation laws.
\newblock {\em IMA J. Numerical Analysis}, 20(1):47--79, 2000.

\bibitem{ToroTitarev2006}
E.F. Toro and V.A. Titarev.
\newblock {MUSTA} fluxes for systems of conservation laws.
\newblock {\em J. Comput. Phys.}, 216(2):403--429, 2006.

\bibitem{TorrilhonExactRS}
M.~Torrilhon.
\newblock Exact {R}iemann solver for ideal {MHD},
  https://web.mathcces.rwth-aachen.de/mhdsolver/.

\bibitem{Torrilhon2002ETH}
M.~Torrilhon.
\newblock Exact solver and uniqueness conditions for {R}iemann problems of
  ideal magnetohydrodynamics.
\newblock Eidgen{\"o}ssische Technische Hochschule [ETH] Z{\"u}rich. Seminar
  f{\"u}r Angewandte Mathematik, 2002.

\bibitem{Torrilhon2003Convergence}
M.~Torrilhon.
\newblock Non-uniform convergence of finite volume schemes for {R}iemann
  problems of ideal magnetohydrodynamics.
\newblock {\em J. Comput. Phys.}, 192(1):73--94, 2003.

\bibitem{Torrilhon2003Uniqueness}
M.~Torrilhon.
\newblock Uniqueness conditions for {R}iemann problems of ideal
  magnetohydrodynamics.
\newblock {\em J. Plasma Physics}, 69(03):253--276, 2003.

\bibitem{Torrilhon2006}
M.~Torrilhon.
\newblock Two-dimensional bulk microflow simulations based on regularized
  {G}rad's 13--moment equations.
\newblock {\em Multiscale Model. Simul.}, 5(3):695--728, 2006.

\bibitem{Torrilhon2012}
M.~Torrilhon.
\newblock Krylov--{R}iemann solver for large hyperbolic systems of conservation
  laws.
\newblock {\em SIAM J. Sci. Comput.}, 34(4):A2072--A2091, 2012.

\bibitem{VanLeer1979}
B.~van Leer.
\newblock Towards the ultimate conservative difference scheme {V}. {A}
  second-order sequel to {G}odunov's method.
\newblock {\em J. Comput. Phys.}, 32:101--136, 1979.

\end{thebibliography}
%=============================================================================================
\smallskip

\end{document}